\documentclass[%
 aip,cha,
 amsmath,amssymb,
 reprint,%
floatfix
]{revtex4-1}
\pdfoutput=1

\usepackage{graphicx}
\usepackage{dcolumn}
\usepackage{bm}

\usepackage[utf8]{inputenc}
\usepackage[T1]{fontenc}
\usepackage{mathptmx}
\usepackage{etoolbox}

\usepackage[mathscr]{euscript}
\newtheorem{Proposition}{Proposition}

\newtheorem{Remark}{Remark}

\makeatletter
\def\@email#1#2{%
 \endgroup
 \patchcmd{\titleblock@produce}
  {\frontmatter@RRAPformat}
  {\frontmatter@RRAPformat{\produce@RRAP{*#1\href{mailto:#2}{#2}}}\frontmatter@RRAPformat}
  {}{}
}%
\makeatother
\begin{document}

\preprint{AIP/123-QED}

\title[Understanding the role of B-cells in CAR T-cell therapy in leukemia through a mathematical model]{Understanding the role of B-cells in CAR T-cell therapy in leukemia through a mathematical model}
\author{Sergio Serrano}
\author{Roberto Barrio}
\affiliation{IUMA, CoDy and Department of Applied Mathematics, Universidad de Zaragoza, Zaragoza, Spain}

\author{Álvaro Martínez-Rubio}
\affiliation{Department of Mathematics, Universidad de C\'{a}diz, Puerto Real, C\'{a}diz, Spain}
\affiliation{Biomedical Research and Innovation Institute of Cádiz (INiBICA), Hospital Universitario Puerta del Mar, Cádiz, Spain}
\email{alvaro.martinezrubio@uca.es}

\author{Juan Belmonte-Beitia}
\author{Víctor M. Pérez-García}
\affiliation{Mathematical Oncology Laboratory (MOLAB), Departament of Mathematics, Instituto de Matemática Aplicada a la Ciencia y la Ingeniería,
Universidad de Castilla-La Mancha, Ciudad Real, Spain.}

\date{\today}

\begin{abstract}
Chimeric Antigen Receptor T (CAR-T) cell therapy has been proven to be successful against different leukaemias and lymphomas.  This paper makes an analytical  and numerical  study of a mathematical model describing the competition of CAR-T, leukaemias tumor and B cells. Considering its significance in sustaining anti-CD19 CAR T-cell stimulation, we  integrate a B-cell source term into the model. Through stability and bifurcation analyses, we reveal the potential for tumor eradication contingent on the continuous influx of B-cells, uncovering a transcritical bifurcation at a critical B-cell input. Additionally, we identify an almost heteroclinic cycle between equilibrium points, providing a theoretical basis for understanding disease relapse. Analyzing the oscillatory behavior of the system, we approximate the time-dependent dynamics of CAR T-cells and leukemic cells, shedding light on the impact of initial tumor burden on therapeutic outcomes. In conclusion, our study provides insights into CAR T-cell therapy dynamics for acute lymphoblastic leukemias, offering a theoretical foundation for clinical observations and suggesting avenues for future immunotherapy modeling research. 
\end{abstract}

\maketitle

\begin{quotation}
CAR T-cell therapy has quickly become the spearhead of cancer immunotherapy, with unprecedented success in the treatment of hematological malignancies. Contrary to classical chemotherapy, the key element of the therapy is a living cell. Understanding its interactions with tumor cells and its microenvironment is crucial for predicting the response. In this article we consider a mathematical model that describes these interactions and focus on the role of B-cells, which theoretically serve as a endogenous vaccine that keeps CAR T-cells activated. We study the properties of the system, characterize the influence of model parameters on the outcome and highlight the therapeutic potential of acting on the processes of B-cell generation. 
\end{quotation}

\section{Introduction}
\label{intro}

Cancer immunotherapy uses components of the patient's immune system to attack cancer cells selectively, and in recent years, it has come to play an important role in treating some types of cancer \cite{Koury18}. Among all of them, the one that has gained more relevance is the treatment with CAR T (Chimeric Antigen Receptor T) cells \cite{DErrico17}. This treatment consists in extracting T-cells from the patient’s blood, adding an antigen receptor to them in the laboratory and subsequently re-infusing them back into the patient’s body. This modification allows T-cells to recognise this antigen in the tumour cells and to kill them \cite{Feins19}. It has been particularly successful against certain types of leukaemia and lymphoma \cite{Maude18,Schuster19}, especially those whose main target is the CD19 antigen, expressed in virtually all B-cells \cite{Sadelain15}.  

In parallel to the clinical and biomedical advances in this therapy, research has also extended to the mathematical modelling of its dynamics. Mathematical models of CAR T-cell therapy have been appearing since the publication of the first clinical trials, building on previous research on T-cell dynamics \cite{Chaudhury20,Nukala21}.  Most of the published works focus on anti-CD19 CAR and specifically applied to lymphoblastic leukaemia \cite{Mostolizadeh18,Stein19,Khailov20,Barros20,Liu21,Barros21,Valle21,Leon21,Perez21,Martinez21,Derippe22}, while others focus on Non-Hodgkin Lymphoma \cite{Liu21,Mueller21,Owens21,Kimmel21}. These models are mainly ODE-based, including from two to up to twenty equations commonly representing therapy, disease and the healthy counterpart. Few of them make use of actual clinical, longitudinal data \cite{Stein19,Liu21,Mueller21}, so most of the research is yet at the proof-of-concept stage. These studies typically consist of the presentation of the problem and the proposal of a model, followed by numerical simulations or fits when data is available. Only a subset of these works focuses on the mathematical analysis of the system of equations, going further than a discussion of the steady states \cite{Mostolizadeh18,Khailov20,Valle21,Marek}. 

An advantage of interdisciplinary applications of mathematics is the appearance of new, interesting problems from the analytical point of view \cite{Semenova}. It is important for the field to not only propose new models and ideas but to support the progress by means of this more theoretical work \cite{Bi,Dickman}. This uncovers model properties, common assumptions and consequences underlying the diversity of models and provides more elements to the discussion of the optimal way of modelling this system. 

In this work, we continue the study done in \cite{Leon21}, where a combination of analytical study and numerical simulations was done to gain clinical insights from the mathematical model. An important simplification was made in that work, i.e., the flow of B cells from the bone marrow was neglected. It was done because in the first weeks after CAR T injection, the main contribution to the dynamics is the expansion of these cells and their effect on the healthy B and leukaemic cells. But when the flow of B cells from the bone marrow can not be neglected, we should consider the contribution of this flow and in this case, the role of B cells is more significant, since a regular supply of B cells is always happening. Thus, in this work, we complete the mentioned gap by making a detailed analytical and numerical study of the model incorporating a constant B-cell input. We will show that this input can be key to obtaining a successful response and sustaining CAR T-cell activation. We also provide a theoretical conceptualization of some features of the therapy such as the time to relapse, the oscillating regime and the maximum values reached by the tumor.  

The structure of the paper is as follows: Firstly, Section II presents the 
mathematical model under consideration. In Section III we perform a mathematical analysis of the model including the existence, uniqueness and stability of steady states. We obtain bifurcations and study the different regions in parameter space and the behaviour of the system in those regions. Section IV presents a study of the importance of the maximum values (or ``peaks'') of leukaemia cells in the treatment, aiming to identify the main relevant parameters that influence this magnitude. In Section V, a sensitivity analysis of the model was performed. Finally, Section VI discusses the findings and summarises the conclusions.

\section{Model presentation}

A mathematical model for the interaction of CAR T-cells, leukemic cells and B-cells was introduced in \cite{Leon21}. It included three compartments for B-cell development and additional equations for cytokines and neurotoxines. A reduced version of this model focusing on early dynamics was further explored. This reduced version reads as follows:

\begin{subequations}
\label{sys0}
\begin{align}
\label{sys0C} \frac{dC}{dt}&=\rho_C(L+B)C-\frac{1}{\tau_c}C,\\
\label{sys0L} \frac{dL}{dt}&=\rho_LL-\alpha CL,\\
\label{sys0B} \frac{dB}{dt}&=-B\left(\alpha C+\frac{1}{\tau_B}\right).
\end{align}
\end{subequations}

Leukemic cells $L(t)$ have a net proliferation rate $\rho_L > 0$ and die due to encounters with CAR T cells. Parameter $\alpha$ measures the probability (per unit time and cell) of an encounter between CAR T and CD19$+$ cells. The equation for CAR-T cells $C(t)$ involves two proliferation terms of CAR T cells due to stimulation by encounters with their target cells: either $L(t)$ or $B(t)$. Parameter $\rho_C$ measures the stimulation to mitosis after encounters with CD19$+$ cells disseminated throughout the whole body (mostly in the circulatory system). The last term in this equation describes the decay of CAR T cells with a mean lifetime $\tau_C$. 
Finally, the equation for B-cells $B(t)$ has two terms. Since B-cells also express the CD19$+$ antigen, they will be targeted by the CAR T cells, thus the term $\alpha BC$. They have a mean lifetime $\tau_B$, which is present in the last term of this equation.

This approximation omits the input of B-cells from the bone marrow, which can be included to yield the following model:

\begin{subequations}
\label{sys}
\begin{align}
\label{sysC}\frac{dC}{dt}&=C\left(\rho_C(L+B+I_0)-\frac{1}{\tau_c}\right),\\
\label{sysL}\frac{dL}{dt}&=\rho_LL-\alpha CL,\\
\label{sysB}\frac{dB}{dt}&=\frac{I_0}{\tau_I}-B\left(\alpha C+\frac{1}{\tau_B}\right).
\end{align}
\end{subequations}

\begin{figure*}[htb]
\centering
\includegraphics[width=0.8\textwidth]{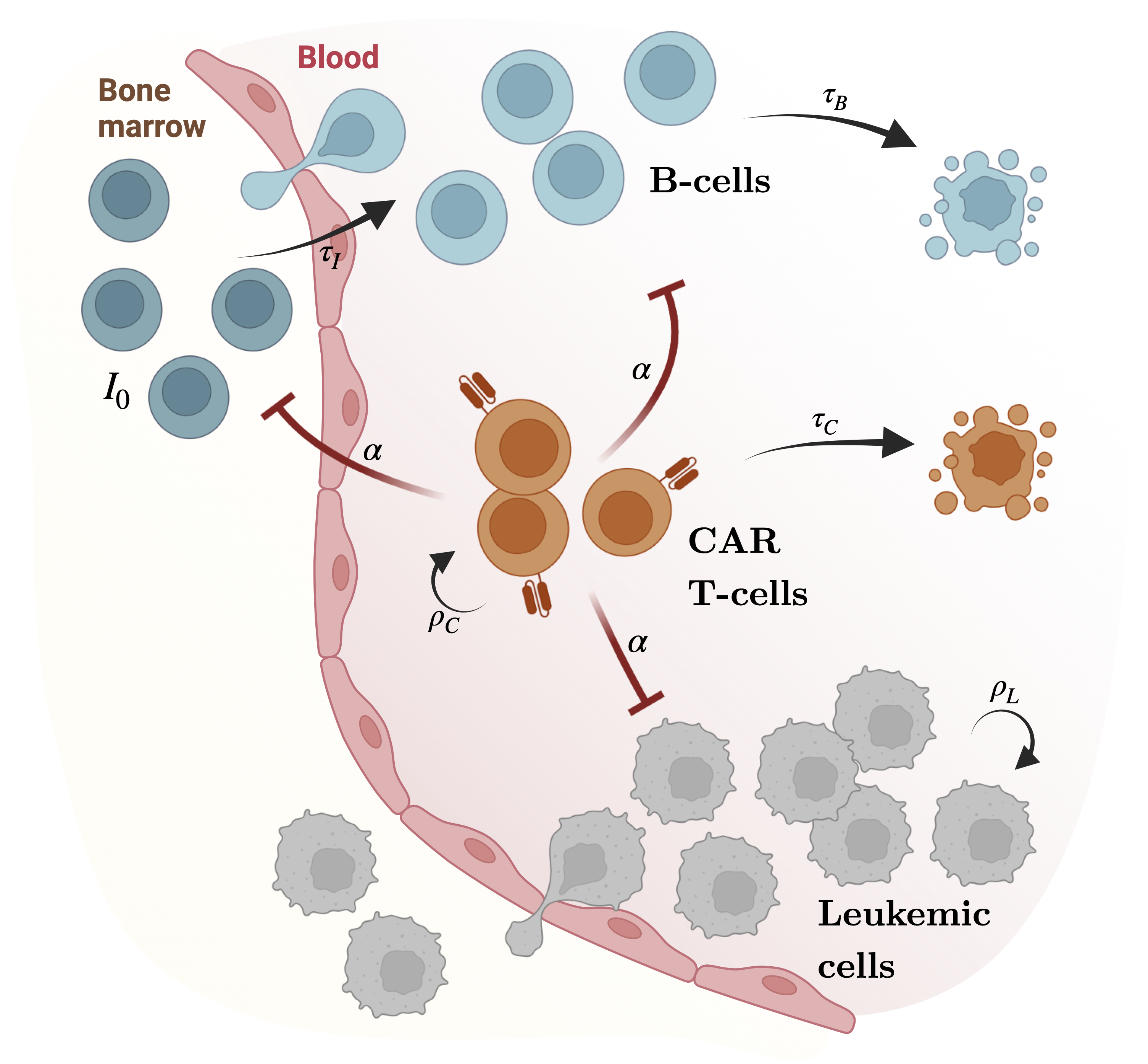}
\caption{Summary of model components and processes. B-cells $B(t)$ exit bone marrow at a rate $I_0/\tau_I$ and die at a rate $1/\tau_B$. Leukemic cells $L(t)$ proliferate with constant rate $\rho_L$. CAR T-cells $C(t)$ proliferate with rate $\rho_C$ upon encounter with antigen-expressing cells, which are removed at a rate $\alpha$. CAR T-cells eventually die at rate $1/\tau_C$.}
\label{diagram}
\end{figure*}

This term is derived from a steady state approximation of a larger version of the model \cite{Leon21}. Note its contribution to the CAR T-cell activation term in the first equation. Parameter $\tau_I$ represents the mean life of B-cell progenitors. Parameter $I_0$ represents the total number of B-cells produced during the development of a B-cell. These cells are also targeted by CAR T-cells, hence its inclusion in the respective equation. Figure \ref{diagram} displays the main processes and elements of the mathematical model.

\subsection{Parameter Selection}

Model \eqref{sys} has several parameters to be estimated. Standard  values are included in Table \ref{tabla}. We describe them briefly here:

\begin{itemize}

\item \textbf{CAR T-cell parameters.} Naive CD8$+$ T cells are quiescent, their mean lifetime ranges from months to years, and they enter the cell cycle following interaction with their antigen \cite{Kasakovski18}. These activated CD8$+$ T cells induce cytolysis of the target cells and secrete cytokines such as TNF$-\alpha$ and IFN$\gamma$. Following activation, most effector cells undergo apoptosis after two weeks, with a small proportion of cells surviving to become CD8$+$ memory T cells capable of longer survival \cite{Kasakovski18}. Recent studies have reported longer survival values of about one month \cite{Ghorashian19}. Thus, we  take the mean lifetime $\tau_C$ of CAR T cells to be in the range of 2-4 weeks. The interaction parameter $\alpha$ and the proliferation rate $\rho_C$ control the timescale and degree of expansion of the CAR T-cell population. Their ranges were determined by simulating behaviours similar to the abundant clinical trial data \cite{Lee15} and in line with values reported in other models \cite{Stein19,Leon21}. Finally, the initial dose is standard \cite{Maude18} and we select an average value of $5\cdot10^7$ cells.

\item \textbf{Leukemic cell parameters} Acute lymphoblastic leukaemias are fast-growing cancers with proliferation rates $\rho_L$ of the order of several weeks. Recent works in lymphopoiesis and CAR T-cell treatment \cite{Martinez21,Chulian21} suggested values between $0.5$ and $0.7$ days$^{-1}$. In those studies, unlike here, the growth is modulated by a signaling function, so the effective rate is lower. We therefore select half the previously estimated value. The initial tumor burden is between 1\% and 100\% of the bone marrow cellularity \cite{Martinez21}, so we select the range $10^{10}-10^{11}$.

\item \textbf{B-cell parameters} Studies on B-cell development yield a B-cell progenitor transition time of around a week \cite{Shahaf16,Chulian21}. From these references we can also compute the stable number of cells produced at the bone marrow $I_0$ and the stable overall level of B-cells $B_0$.

\end{itemize}

\begin{table*}
\centering
\begin{ruledtabular}
\begin{tabular*}{\linewidth}{@{\extracolsep{\fill}} lccccc}
Parameter & Description & Standard Value & Range & Units & Source\\ 
\hline
$\rho_C$ & Stimulation of CAR T-cells & $10^{-11}$ & $5\cdot10^{-12}-5\cdot10^{-11}$ & day$^{-1}\cdot$ cell$^{-1}$ &  \cite{Leon21}\\
$\tau_C$ & Activated CAR T-cell lifetime & $20$ & $14-30$ & day & \cite{Ghorashian19} \\
$\rho_L$ & Leukaemic cell growth rate & $0.2$ & $0.1-0.3$ & day$^{-1}$ & \cite{Martinez21,Skipper70}\\
$\alpha$ & Killing efficiency of CAR T-cells & $10^{-11}$ & $5\cdot10^{-12}-5\cdot10^{-11}$   & day$^{-1}\cdot$ cell$^{-1}$  & \cite{Leon21}\\
$I_0$ &  Bone marrow B-cell output & $10^9$ & $5\cdot10^8 - 5\cdot10^9$ & cell & \cite{Chulian21} \\
$\tau_I$ &  B-cell maturation time & $4$ & $1-7$ & day & \cite{Chulian21,Shahaf16}\\
$\tau_B$ & B-cell lifetime & $45$ & $30-60$ & day & \cite{Fulche97}\\
$L_0$ &  Initial tumor burden & $5\cdot10^{10}$ & $10^{10}-10^{11}$ & cell & \cite{Martinez21} \\
$C_0$ &  CAR T-cell dose & $5\cdot10^7$ & $10^7-10^8$ & cell & \cite{Maude18} \\
$B_0$ &  Initial B-cell load & $5\cdot10^8$ & $10^8-10^9$ & cell & \cite{Chulian21} \\
\end{tabular*}
\end{ruledtabular}
\caption{Parameter values for system of equations \eqref{sys}}\label{tabla}
\end{table*}

\section{Model analysis}
\label{sec:analisis}

In this section we make an analysis of some basic properties of system~\eqref{sys}.

First, we state the biological sense of the model, that is, if the model is well-posed in the sense that initial positive values give positive solutions of the system.

\begin{Proposition}
Let the parameters $I_0, \rho_C, \rho_L, \alpha$ be non-negative and $\tau_c, \tau_B, \tau_I$ positive in system \eqref{sys}. Then, for any non negative initial data, the solutions of the system are positive.
\end{Proposition}

\emph{Proof}: 
Just take into account that $L=0$ and $C=0$ are invariant planes of ~\eqref{sys} and on $B=0$ the vector field points towards the first octant.

Regarding the third equation, let $\textbf{n}_3=(0,0,-1)$ denote the outward normal unit vector to plane $B=0$. Consider the scalar product of the ODE system \eqref{sys} with $\textbf{n}_3$ at hyper-surface $B=0$:
$$
\left\langle \left(\frac{dC}{dt},  \frac{dL}{dt},  \frac{dB}{dt}\right), (0,0,-1) \right\rangle_{B=0}=-\frac{I_0}{\eta_I}<0.
$$
Thus, the piece of hyper-surface $\{B=0\}\cap \mathbb{R}^3_{+}$ is semipermeable inward $\mathbb{R}^3_+$. 

As a result, $\mathbb{R}^3_+$ is a positively invariant domain for system \eqref{sys}. Therefore, non-negativity of solutions follows.

\rule{5pt}{5pt}

\subsection{Stability of steady states}

Next, we look for steady states. Equating the previous system to zero we obtain the following equilibrium points:

\begin{subequations}
\label{eqpoints}
\begin{align}
\label{eqpoints1} P_1&=\left(0,0,\frac{I_0\tau_B}{\tau_I}\right),\\
\label{eqpoints2} P_2&=\displaystyle{\left(\frac{1}{\alpha}\left(\frac{I_0}{\tau_I(\frac{1}{\tau_C\rho_C}-I_0)}-\frac{1}{\tau_B}\right),0,\frac{1}{\tau_C\rho_C}-I_0\right)},\\
\label{eqpoints3} P_3&=\displaystyle{\left(\frac{\rho_L}{\alpha},
      \frac{1}{\tau_C\rho_C}-I_0\left(1+\frac{\tau_B}{\tau_I(1+\tau_B\rho_L)}\right),\frac{I_0}{\tau_I\left(\rho_L+\frac{1}{\tau_B}\right)}
\right).  }
\end{align}
\end{subequations}

The first equilibrium point $P_1$ describes a scenario in which CAR T-cells have completely eliminated leukemic cells and they have eventually died out, resulting in a stable, homeostatic level of B-cells in which the bone marrow output $I_0$ compensates the natural death rate of B-cells $\tau_B$. The second equilibrium point $P_2$ again describes a complete cure situation, but in this case CAR T-cells coexist to B-cells due to the stimulation that the latter provide to the former. Lastly, equilibrium point $P_3$ describes a situation in which the three cell types coexist. To study the local stability of these steady states, we calculate Jacobian matrix of Eqs. \eqref{sys}:
\begin{equation}
J(C,L,B)=
\begin{pmatrix}
\rho_C(L+B+I_0)-\frac{1}{\tau_C} & \rho_C C & \rho_C C\\
-\alpha L & \rho_L-\alpha C & 0\\
-\alpha B & 0 & -\alpha C-\frac{1}{\tau_B}
\end{pmatrix}.
\end{equation}
and obtain the eigenvalues. For equilibrium point $P_1$ we find 
\begin{flalign*}
   \lambda_1&=\rho_L>0, &\\
   \lambda_2&=-\frac{1}{\tau_B}<0,&\\
   \lambda_3&=\rho_CI_0\left(1+\frac{\tau_B}{\tau_I}\right)-\frac{1}{\tau_C}.&
\end{flalign*}
Clearly, $P_1$ is a saddle point, and, therefore, unstable. Moreover, if

\begin{equation}\label{tr1}
I_0>\displaystyle{\frac{1}{\tau_C\rho_C}\left(1-\frac{\tau_B}{\tau_B+\tau_I}\right)},
\end{equation}
the dimension of the unstable manifold $W^u(P_1)$ becomes $2$. Otherwise, its dimension is equal to $1$. For equilibrium point $P_2$ we find that all the components are non-negative if and only if 
\begin{equation}
I_0\in \left[\frac{1}{\tau_C\rho_C}\left(1-\frac{\tau_B}{\tau_B+\tau_I}\right),\frac{1}{\tau_C\rho_C}\right).
\end{equation}
The corresponding eigenvalues are

\begin{align*}
\lambda_1 &=\rho_L-\alpha C_2,\\
\lambda_{2,3} &=\frac{1}{2}\left(\frac{-I_0}{\tau_I\left(\frac{1}{\tau_C\rho_C}-I_0\right)}
\pm\sqrt{\frac{I_0^2}{\tau_I^2\left(\frac{1}{\tau_C\rho_C}-I_0\right)^2}-4\alpha\rho_CB_2C_2} \right).
\end{align*}

$P_2$ is asymptotically stable if and only if all eigenvalues have negative real part. Then, on the one hand, if
\begin{eqnarray*}
  \rho_L-\alpha C_2&=&\rho_L+\frac{1}{\tau_B}-\frac{I_0}{\tau_I(\frac{1}{\tau_C\rho_C}-I_0)}<0 \Leftrightarrow\\
  &\Leftrightarrow& I_0\in \left(\frac{1}{\tau_C\rho_C}\left(1-\frac{\tau_B}{\tau_B+\tau_I(1+\tau_B\rho_L)}\right),\frac{1}{\tau_C\rho_C}\right).
\end{eqnarray*}
On the other hand, if $B_2C_2>0$, the real part of the eigenvalues $\lambda_{2,3}$ has the same sign. Moreover,  $\lambda_{2,3}$ are real numbers if 

\begin{align*}
 \frac{4}{\tau_B\tau_C^3\rho_C^2}&-\frac{4I_0}{\tau_C^2\rho_C}\left(\frac{1}{\tau_I}+\frac{3}{\tau_B}\right)+
  I_0^2\left(\frac{1}{\tau_I^2}+\frac{8}{\tau_I\tau_C}+\frac{12}{\tau_B\tau_C}\right)\\
  &-4\rho_CI_0^3\left(\frac{1}{\tau_I}+\frac{1}{\tau_B}\right)>0.
\end{align*}

Therefore, $P_2$ is asymptotically stable if
\begin{equation}
I_0\in \left(\frac{1}{\tau_C\rho_C}\left(1-\frac{\tau_B}{\tau_B+\tau_I(1+\tau_B\rho_L)}\right),\frac{1}{\tau_C\rho_C}\right).
\end{equation}

In the same way than the previous case, all the components of the equilibrium point $P_3$ are non-negative if and only if
\begin{equation}\label{condition_1}
I_0\in \left[ 0,\frac{1}{\tau_C\rho_C}\left(1-\frac{\tau_B}{\tau_B+\tau_I(1+\tau_B\rho_L)}\right)\right].
\end{equation}
We denote this critical value of $I_0$ as
\begin{equation}
\label{transcb}
I^{\text{crit}}_0=\frac{1}{\tau_C\rho_C}\left(1-\frac{\tau_B}{\tau_B+\tau_I(1+\tau_B\rho_L)}\right).
\end{equation}
The eigenvalues result from obtaining the roots of the characteristic polynomial

\begin{align*}
   p(\lambda)&=\lambda^3+\lambda^2\left(\frac{1}{\tau_B}+\rho_L\right)+\lambda\rho_C\rho_L\left(\frac{1}{\rho_C\tau_C}-I_0\right)\\
   &+\rho_C\rho_L L_3\left(\frac{1}{\tau_B}+\rho_L\right).
\end{align*}
    
Using the Routh-Hurwitz criterion, we can state that if condition \eqref{condition_1} is satisfied, then the equilibrium point $P_3$ is asymptotically stable.

\subsection{Bifurcations}

Considering the non-hyperbolic equilibrium points, we detect several bifurcations that we can state as follows:

\begin{Proposition}
Let $P_3$ be the equilibrium point above described,  and consider $I_0=0$. Then, for this value of the parameter, we obtain a \emph{subcritical Hopf bifurcation}, since the first Lyapunov coefficient is $l_1=\tau_C^{3/2}\rho_C^2/(3\sqrt{\rho_L})>0$. Hence, unstable periodic orbits are generated for $I_0>0$.
\end{Proposition}

\begin{Proposition}
Let $P_1$ and $P_2$ be the equilibrium points above described,  and consider $I_0=\displaystyle{\frac{1}{\tau_C\rho_C}\left(1-\frac{\tau_B}{\tau_B+\tau_I}\right)}$ (\ref{tr1}). Then, at this value of the parameter,
$P_1$ and $P_2$ develop a \emph{transcritical bifurcation}: Both equilibrium points exchange the dimensions of their stable and unstable manifolds.
\end{Proposition}

\begin{Proposition}
Let $P_2$ and $P_3$ be the equilibrium points above described,  and consider $I_0=I^{\text{crit}}_0$ (\ref{transcb}). Then, at this value of the parameter,
$P_2$ and $P_3$ develop a \emph{transcritical bifurcation}: $P_2$ becomes a stable point while $P_3$ becomes an unstable point when $I_0>I^{\text{crit}}_0$.
\end{Proposition}

\begin{figure*}[ht!]
\centering
\includegraphics[width=0.9\textwidth]{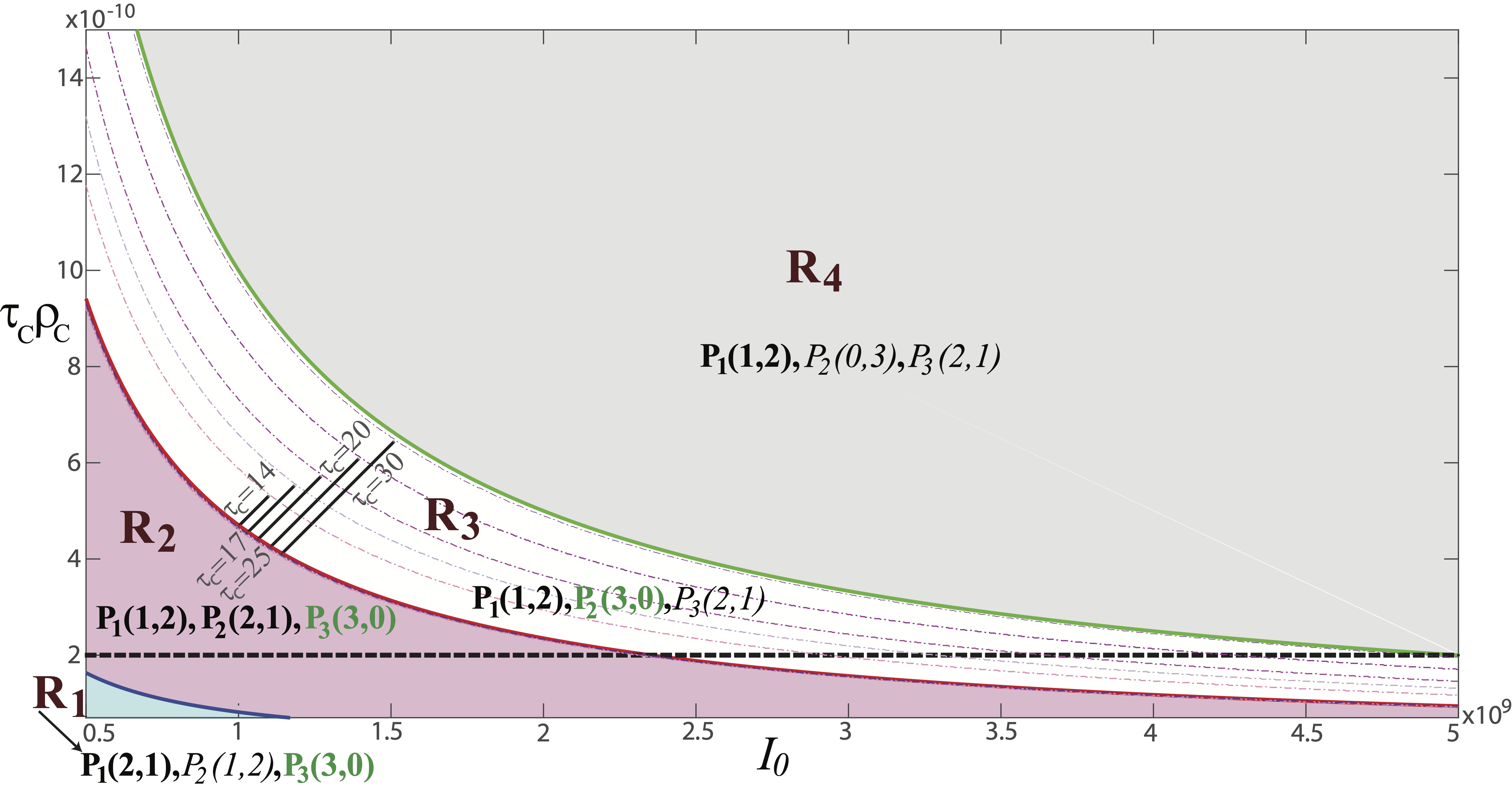}
\caption{Regions determined by the stability of the equilibrium points with $\rho_L=0.2$ day$^{-1}\cdot$ cell$^{-1}$; 
$\alpha=10^{-11}$ day$^{-1}\cdot$ cell$^{-1}$; $\tau_I=4$ days; $\tau_B=45$ days. 
Bold text denotes equilibrium points with biological meaning (all components are positive). Green text denotes stable equilibria. Dashed black line represent the standard value of $\tau_C\rho_C=2\cdot 10^{-10}$ cell$^{-1}$. See main text for more details.}
\label{curvas}
\end{figure*}

In Figure~\ref{curvas} we study in the biparametric plane $(I_0, \, \tau_C\rho_C)$ the basic bifurcations that delimit different regions in function of the stability properties of the equilibrium points. We have located four possible regions ($R_1$ to $R_4$) depending on the different number of equilibria and their topological characteristic. In particular, we represent the following curves (given by the previous analysis):
\begin{itemize}
\item Blue curve: $\displaystyle{I_0=\frac{1}{\tau_C\rho_C}\left(1-\frac{\tau_B}{\tau_B+\tau_I}\right)}$;
  \item Red curve: $\displaystyle{I_0=\frac{1}{\tau_C\rho_C}\left(1-\frac{\tau_B}{\tau_B+\tau_I(1+\tau_B\rho_L)}\right)}$;
  \item Green curve: $\displaystyle{I_0=\frac{1}{\tau_C\rho_C}}$;
  \item Finally, we represent the standard value of  $\tau_C\rho_C=2\cdot 10^{-10}$, in a dashed horizontal black curve.
\end{itemize}

Analyzing these regions we obtain some interesting results for the model, which we summarize as follows:
\begin{itemize}
  \item Region $R_1$:
    \begin{itemize}
      \item $P_1$ has biological sense (all the coordinates are positive) and it is an unstable point. The stable manifold $W^s(P_1)$ has dimension $2$ ($\textrm{dim}(W^s(P_1))=2$) and the unstable manifold $W^u(P_1)$ has dimension $1$ ($\textrm{dim}(W^u(P_1))=1$). So, its characteristic dimensions are $(\textrm{dim}(W^s(P_1)), \, \textrm{dim}(W^u(P_1))) = (2,1)$.
      \item $P_2$ does not have biological sense (at least one coordinate is negative) and it is an unstable point. Now  $(\textrm{dim}(W^s(P_2)), \, \textrm{dim}(W^u(P_2))) = (1,2)$.
      \item $P_3$ has biological sense and it is an asymptotically stable equilibrium point. So, now  $(\textrm{dim}(W^s(P_3)), \, \textrm{dim}(W^u(P_3))) = (3,0)$.
    \end{itemize}
  \item Region $R_2$:
    \begin{itemize}
      \item $P_1$ has biological sense and it is an unstable point. Now $(\textrm{dim}(W^s(P_1)), \, \textrm{dim}(W^u(P_1))) = (1, 2)$.
      \item $P_2$  has biological sense and it is an unstable point. Now   $(\textrm{dim}(W^s(P_2)), \, \textrm{dim}(W^u(P_2))) = (2,1)$.
      \item $P_3$ has biological sense and it is an asymptotically stable equilibrium point. So, now $(\textrm{dim}(W^s(P_3)), \, \textrm{dim}(W^u(P_3))) = (3, 0)$.
    \end{itemize}
  \item Region $R_3$:
    \begin{itemize}
      \item $P_1$ has biological sense and it is an unstable point. Now $(\textrm{dim}(W^s(P_1)), \, \textrm{dim}(W^u(P_1))) = (1, 2)$.
      \item $P_2$ has biological sense and it is an asymptotically stable equilibrium point, so, now  $(\textrm{dim}(W^s(P_2)), \, \textrm{dim}(W^u(P_2))) = (3, 0)$.
      \item $P_3$ does not have biological sense and it is an unstable point. Now $(\textrm{dim}(W^s(P_3)), \, \textrm{dim}(W^u(P_3))) = (2,1)$.
    \end{itemize}
  \item Region $R_4$:
    \begin{itemize}
      \item $P_1$  has biological sense and it is an unstable point. Now  $(\textrm{dim}(W^s(P_1)), \, \textrm{dim}(W^u(P_1))) = (1,2)$.
      \item $P_2$ does not have biological sense and it is an unstable point. Now  $(\textrm{dim}(W^s(P_2)), \, \textrm{dim}(W^u(P_2))) = (0, 3)$.
      \item $P_3$ does not have biological sense and it is an unstable point. Now  $(\textrm{dim}(W^s(P_3)), \, \textrm{dim}(W^u(P_3))) = (2, 1)$.
          \item Note that in this region $\rho_CI_0-\frac{1}{\tau_C}$ is positive. Therefore,  if we start with positive initial conditions, $\frac{{\rm d}C}{{\rm d}t}\geq k\, C$ with $k$ a positive constant. So $C$ grows indefinitely to infinity. This, in turn, causes $L$ to converge to $0$. Thus, for this case, the solutions are unbounded.
    \end{itemize}
\end{itemize}
Note that in some regions only some equilibrium points have biological sense, whereas others have at least one negative coordinate. In short, $P_1$ is always biologically meaningful; $P_2$ is only acceptable in regions $R_2$ and $R_3$, and $P_3$ is only acceptable in regions $R_1$ and $R_2$. This means region $R_2$ is the only one in which the three equilibrium points are available, which makes it more interesting from the biological point of view. Stability also differs across regions; $P_3$ is the stable point in regions $R_1$ and $R_2$, until $I_0=I_0^{crit}$, when $P_3$ loses biological sense and $P_2$ becomes the stable point (in region $R_3$). Finally, in region $R_4$, all equilibrium points are unstable. 

\begin{figure*}[!]
\centering
\includegraphics[width=0.9\textwidth]{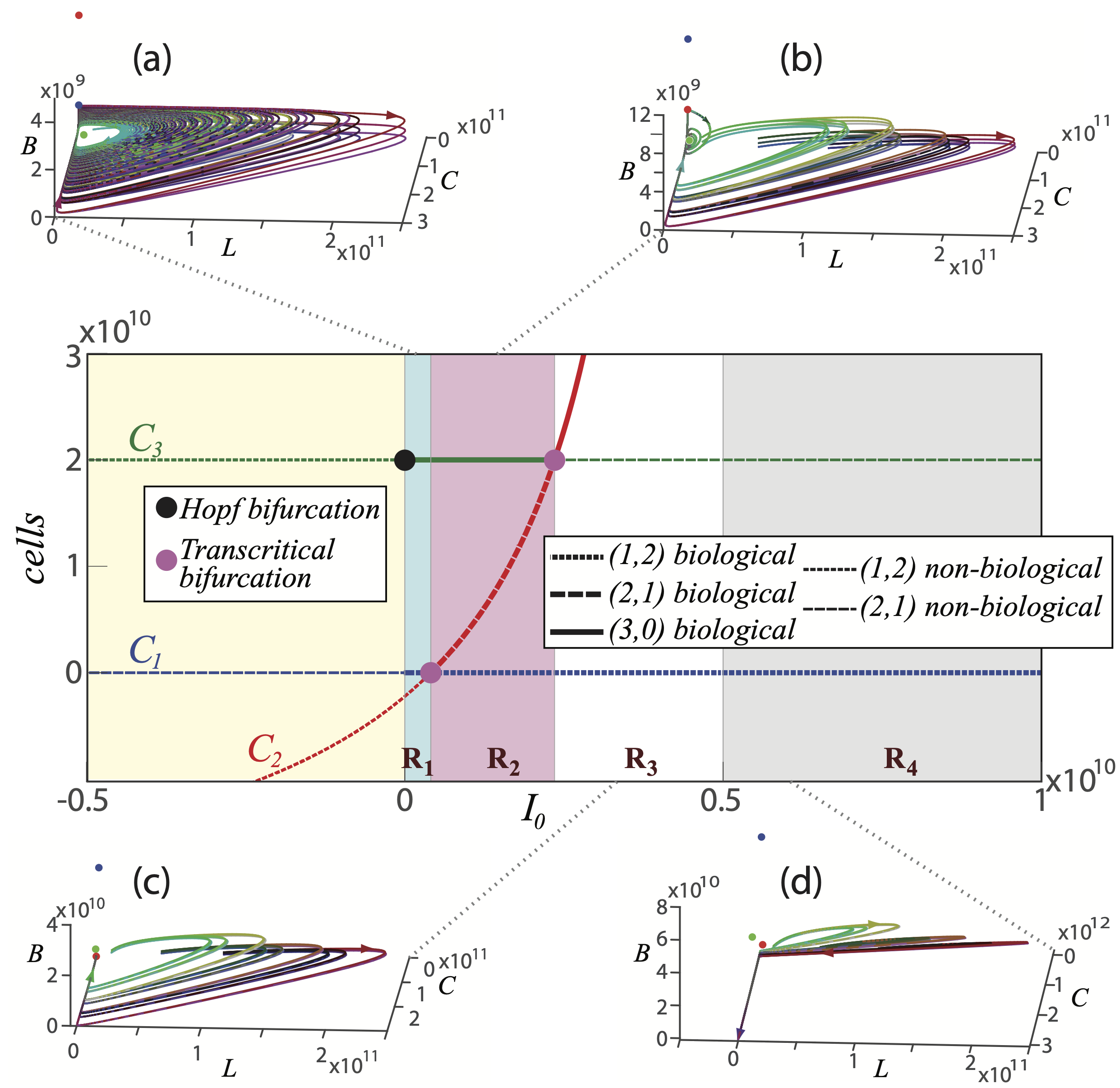}
\caption{Central panel, value of the components $C_i$, $i=1 \ldots 3$, of the equilibrium points as a function of the B-cell influx $I_0$, for standard values of the parameters ($\rho_C\tau_C=2\cdot 10^{-10}$ cell$^{-1}$; $\rho_L=0.2$ day$^{-1}\cdot$ cell$^{-1}$, $\alpha=10^{-11}$ day$^{-1}\cdot$ cell$^{-1}$, $\tau_I=4$ days, $\tau_B=45$ days). 
The interval for $I_0$ is broader than the biological one, to incorporate the bifurcations happening at $I_0=0$. Dashed lines indicate that the equilibrium is unstable, while solid lines indicate that the equilibrium is stable. Also, bold lines denote that the equilibrium point is biologically meaningful. The black circle marks the Hopf bifurcation of $P_3$. The two purple circles mark transcritical bifurcations. Panels (a) to (d) show a characteristic phase diagram for several close initial conditions in
 the different regions of equilibrium ($R_1$ to $R_4$ respectively).}
\label{equilibrios}
\end{figure*}

To further analyse the dynamics of the system, Figure \ref{equilibrios} shows the standard behaviour that occurs in each region. This behaviour was calculated for standard parameter values (Table \ref{tabla}), in particular $\rho_C \tau_C=2\cdot 10^{-10}$ cell$^{-1}$ (black dashed line in Figure \ref{curvas}), but the qualitative behaviour is the same over the biological range of parameters. Note that on the left hand side (with negative $I_0$) none of the equilibrium points make biological sense and all are unstable. 
In $R_1$, the trajectory converges to $P_3$, which is the only stable equilibrium, but it evolves by first passing very close to $P_1$. In $R_2$ the situation is similar, but this time it passes very close to $P_2$. In $R_3$ it converges to $P_2$, the only stable equilibrium point. And in the $R_4$ region, the evolution is not bounded in variable $C$. Note that in all cases, the first thing that is observed is a strong peak in the number of tumour cells. 

\begin{Remark}
If the initial condition $C_0$ is zero, then $C(t)\equiv 0$; $L(t)=e^{\rho_L t}L_0$ and $B(t)=\frac{I_0\tau_B}{\tau_I}+e^{-\frac{t}{\tau_B}}\left(B_0-\frac{I_0\tau_B}{\tau_I}\right)$.
Thus, $L$ grows indefinitely and $B$ converges to $\frac{I_0\tau_B}{\tau_I}$. This behaviour occurs for any regions determined in Fig. \ref{curvas}. The biological interpretation is that, in absence of therapy, the tumor cell population would grow indefinitely, up to a point which would be fatal for the patient. 
\end{Remark}

\begin{figure*}[!]
\centering
\includegraphics[width=0.9\textwidth]{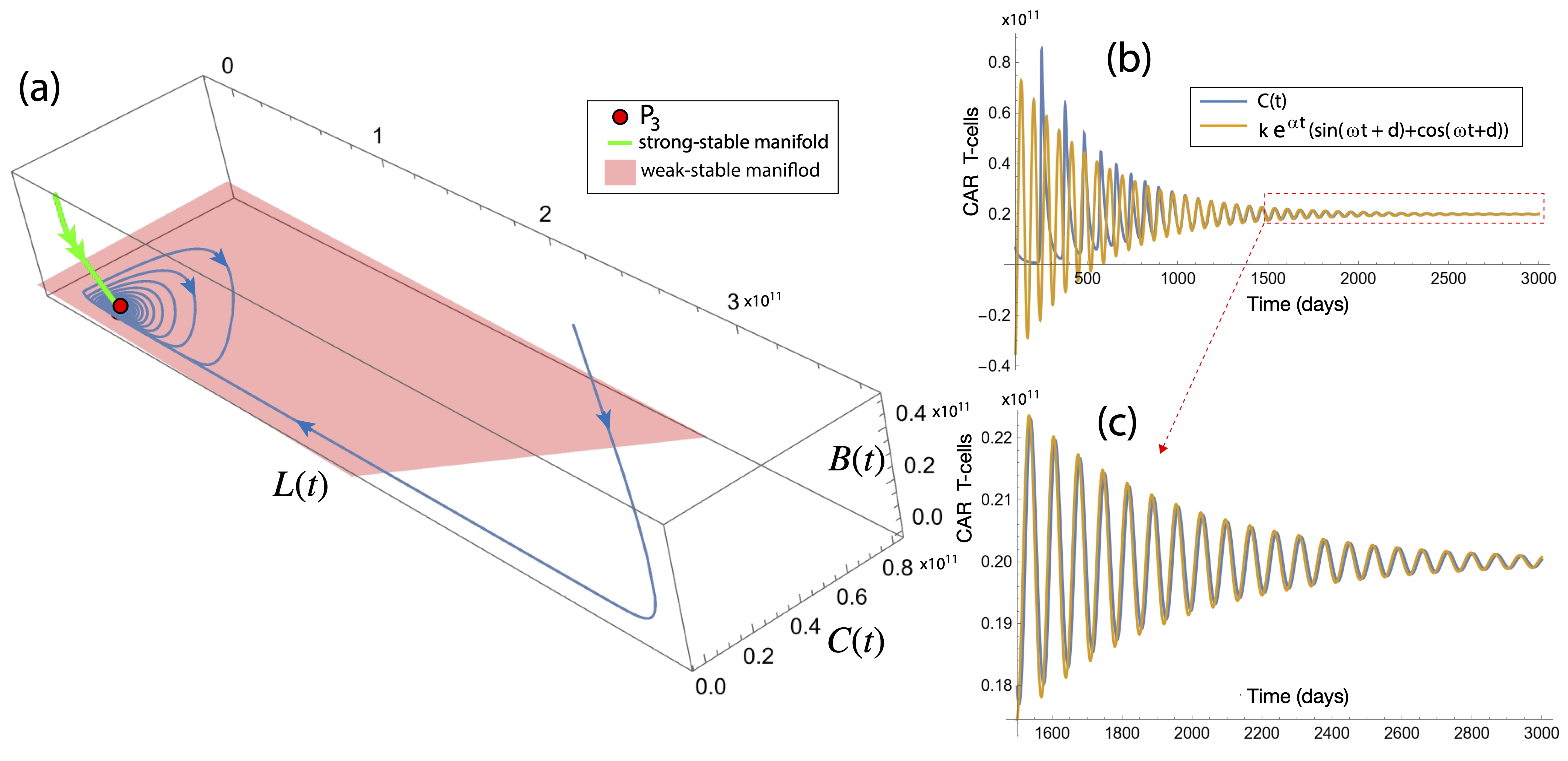}
\caption{Simulations from Eqs. \eqref{sys} displaying focus-like dynamics in Region $R_2$ approaching $P_3$ equilibrium. Parameter values for the orbit shown (blue curve) are $I_0= 0.5 \cdot 10^9$ cells, $\rho_C \tau_C = 2 \cdot 10^{-10}$ cell$^{-1}$, $\rho_L=0.2$ day$^{-1}\cdot$ cell$^{-1}$, 
$\alpha=10^{-11}$ day$^{-1}\cdot$ cell$^{-1}$, $\tau_I=4$ days, $\tau_B=45$ days. The initial state is given by $B_0=5\cdot10^8$ cells; $L_0=5\cdot10^{10}$ cells and $C_0=5\cdot10^7$ cells.  \textbf{(a)} Phase-portrait around $P_3$ highlighting strong and weak stable manifolds. \textbf{(b,c)} Evolution of CAR T-cells with time asymptotically approximated by $x(t)$ (Eqn.~(\ref{approxeq})) with $K=1/2 \cdot 10^{11}, k_s=k_c=1, d=4.8$.}
\label{focoDef}
\end{figure*}

Going back to Figure \ref{curvas}, besides region-delimiting curves we have also included families of curves that take into account changes in the invariant manifolds of $P_3$ equilibrium: We select a series of values for $\tau_C$ ($14, 17, 20, 25$ and $30$ days) and for each value two discontinuous curves are plotted. Note that the left limit curves are all very close each other. For values between these two curves, all eigenvalues of $P_3$ are real. A diagonal segment connecting both curves has been drawn for each $\tau_C$ to make clearer to the reader the corresponding limit curves. Outside this range, $P_3$ has two complex conjugate eigenvalues. We can see how the entire $R_1$ region and almost the entire $R_2$ region are within the lower zone where $P_3$ has a stable focus-node behavior. The strong stable manifold is given by the real negative eigenvalue, and therefore the orbit approaches fast the
2D weak stable manifold where a typical stable focus dynamics is developed. Thus, once in this manifold the dynamics can be very well approximated by the formula
\begin{equation}
\label{approxeq}
x(t) = K \cdot e^{\alpha t} (k_s \sin(\omega t + d) + k_c \cos(\omega t + d))+x_3,
\end{equation}
where $\alpha = \Re(\lambda_{1,2})$ and $\omega  = |\Im(\lambda_{1,2})|$, are the real part and absolute value of imaginary part of the complex eigenvalues, respectively. $x_3$ is the corresponding component of $P_3$ and $K, k_s, k_c, d  \in \mathbb{R}$. In Figure~\ref{focoDef} we illustrate such behaviour with one orbit
for $I_0= 0.5 \cdot 10^9$ cells and using $\rho_C \tau_C = 2 \cdot 10^{-10}$ cell\textsuperscript{-1}. The orbit oscillates from the beginning and it takes some time to be well approximated by the focus dynamics (the time to go to the weak stable manifold), but, once inside it, it is well approximated by (\ref{approxeq}) using $K=1/2 \cdot 10^{11}, k_s=k_c=1, d=4.8$.

\begin{figure*}[!]
\centering
\includegraphics[width=0.9\textwidth]{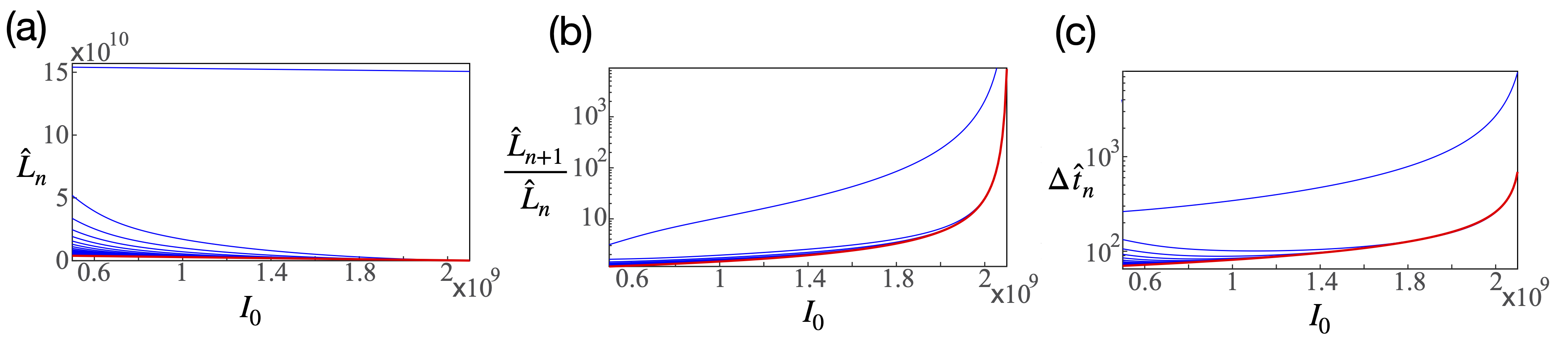}
\caption{Asymptotic behaviour of leukaemic cells approaching equilibrium point $P_3$ as a function of $I_0$, for the same parameter values as in Fig. \ref{equilibrios}. \textbf{(a)} Magnitude of the local maxima of $L(t)$ approaching $L_3$ (red curve). \textbf{(b)} Quotient between two consecutive local maxima of $L$. Red curve marks $e^{-(2\pi )\alpha/\omega}$. \textbf{(c)} Difference between the times at which two consecutive local maxima of $L$ occur. Red curve marks $(2\pi )/\omega$. In all three cases, the curves appear in descending order with time. \label{fococur}}
\end{figure*}

\begin{figure*}[ht!]
\centering
\includegraphics[width=0.9\textwidth]{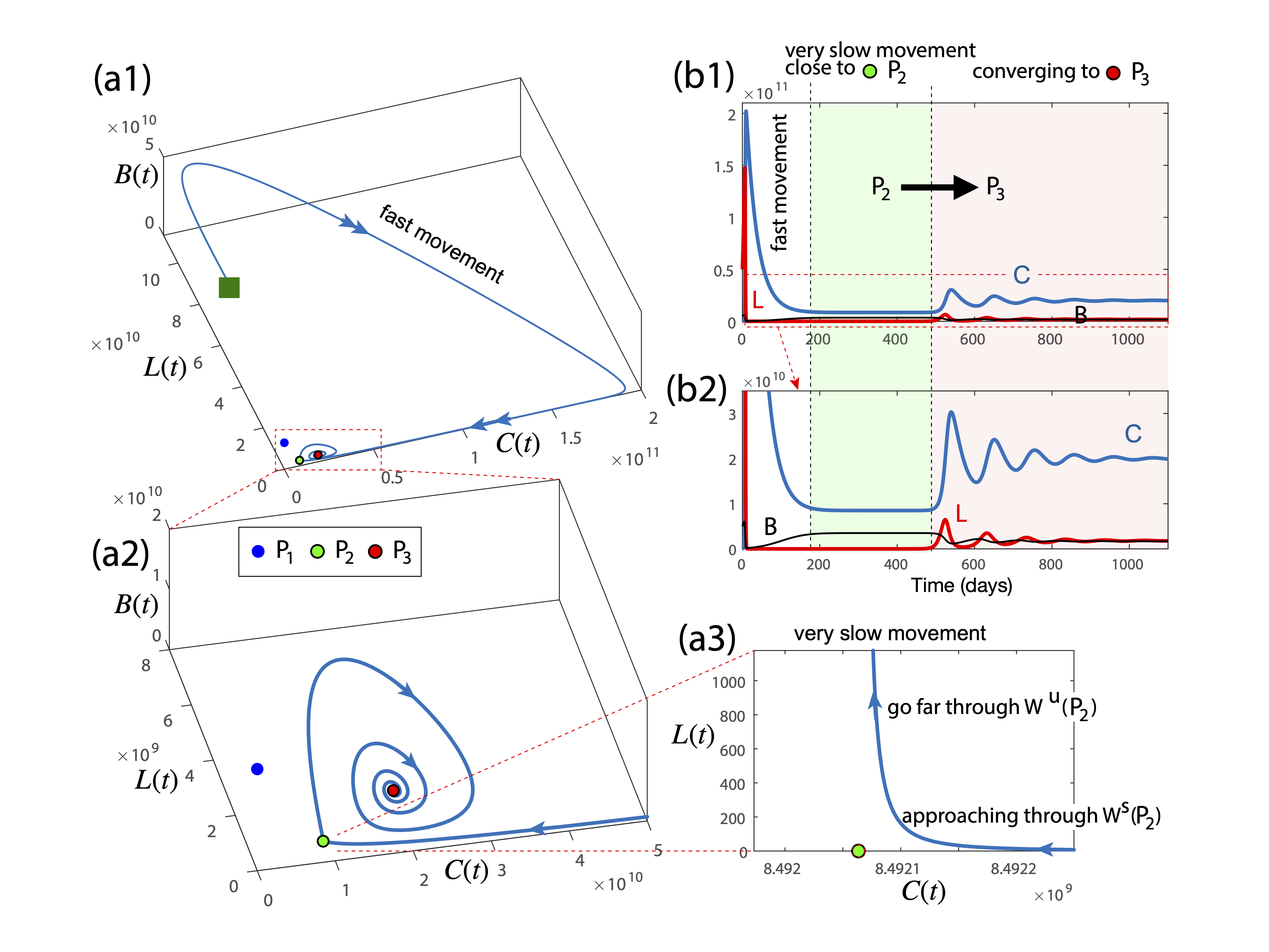}
\caption{Simulation of system \eqref{sys} showing an almost heteroclinic cycle between $P_2$ and $P_3$ equilibria in the $R_2$ region close to the $R_3$ region. Parameter values are $\rho_C \tau_C = 2 \cdot 10^{-10}$ cell$^{-1}$, $I_0=1.5\cdot10^9$ cells, $\rho_L=0.2$ day$^{-1}\cdot$ cell$^{-1}$, $\alpha=10^{-11}$ day$^{-1}\cdot$ cell$^{-1}$, $\tau_I=4$ days, $\tau_B=45$ days. \textbf{(a)} Phase portrait including all three equilibria of the system. Green square denotes initial conditions. Subfigures (a2) and (a3) represent magnifications of the phase plane. \textbf{(b)} Evolution in time of the three variables of the system showing the transition from $P_2$ to $P_3$.}
\label{het2eq}
\end{figure*}

\begin{figure*}[!]
\centering
\includegraphics[width=0.9\textwidth]{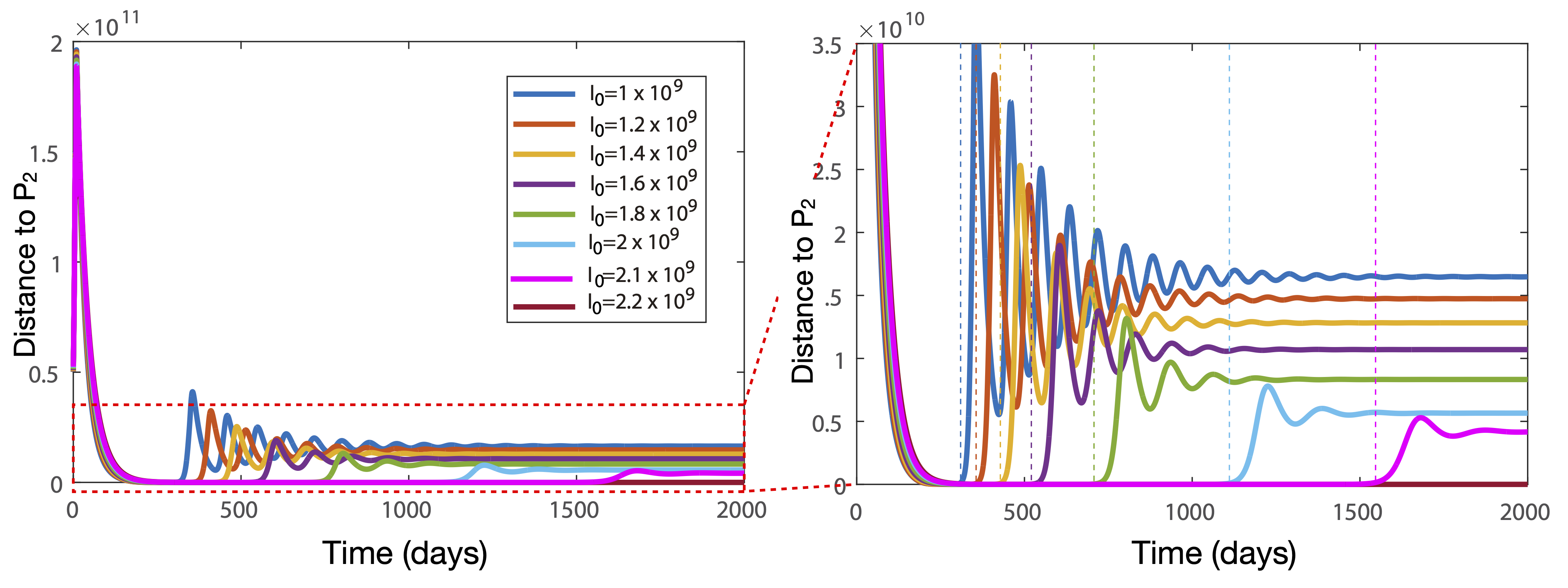}
\caption{Distances of the orbits to the equilibrium $P_2$ for different values of $I_0$ inside the $R_2$ region. The rest of parameters are as in Figure \ref{het2eq}.}
\label{distancia}
\end{figure*}

In Figure~\ref{fococur} we can see that the evolution towards focus dynamics is general in this region once enough time has elapsed to fall into the 2D weak manifold. In this figure, $I_0$ is left as a free parameter and the other parameters take the standard values. The left panel represents the values of the different local maxima of the variable $L$ that the orbit reaches over time. We denote this maxima as $\hat{L}_n$. The upper curve is the first local maximum, the following blue curves represent the following local maxima, in decreasing order towards the value of $L_3$ represented by the red curve. The central panel shows the ratio of successive maxima. The theoretical value given by the Eq.~(\ref{approxeq}) is $x(t)/x(t+2\pi/\omega)=e^{-(2\pi )\alpha/\omega}$. This value is represented by the red curve. For all values of $I_0$ the ratio converges to the theoretical value. Finally, the right panel shows the difference between the times at which two consecutive local maxima of $L$ occur, denoted by $\Delta \hat{t}_n$ The theoretical value given by the Eq.~(\ref{approxeq}) is $(2\pi )/\omega$, represented by the red curve. Again, we can see how the actual values monotonically approach the theoretical value.

We remarked that region $R_2$ is the only one with three biologically meaningful equilibrium points, so the dynamics are expected to be richer. These are mainly determined by the interplay between the manifolds of the $P_2$ and $P_3$ equilibria. Near the upper limit of the $R_2$ region there is an almost  heteroclinic connection among both equilibria, giving rise to a dynamics dominated initially by the $P_2$ equilibrium, approaching it by its 2D stable manifold and leaving by its 1D unstable manifold. Once it is close enough to $P_3$ it goes fast by the strong stable manifold to the 2D weak stable manifold and begins to behave like a stable focus dynamics. This behaviour can be seen in Figure~\ref{het2eq}(a). We first have a large increment in all variables to go later close to the $P_2$, remaining there a long time (this means a large time with a low value of leukaemic cells $L$), but later the dynamics goes to the endemic equilibrium point $P_3$ (see the magnifications (a2) and (a3)).
The figures on the right ((b1) and (b2)) show the time series of the solution showing the long time spent close to $P_2$ (the rest of the dynamics is very fast). The last part of the time series illustrates the stable focus dynamics. Therefore, in some situations large oscillations are expected (see the first large oscillation) and large excursions from the neighborhood of one equilibrium to another one.

To explore deeper the relationship between the distance to the equilibrium point $P_2$ and the time spent in its neigborhood, i.e. the time with low $L$, in Figure \ref{distancia} we simulate the dynamics of the system \eqref{sys} for different values of $I_0$ inside the $R_2$ region. We clearly observe how approaching the $R_3$ region the behaviour is very close to $P_2$ for a long time. This means that for a long time interval the illness may not be present at detectable levels, in some cases for more than 5 years (e.g. $I_0=2.1\cdot 10^{9}$ and $2.2 \cdot 10^{9}$ cells).

\begin{figure*}[!]
\centering
\includegraphics[width=0.85\textwidth]{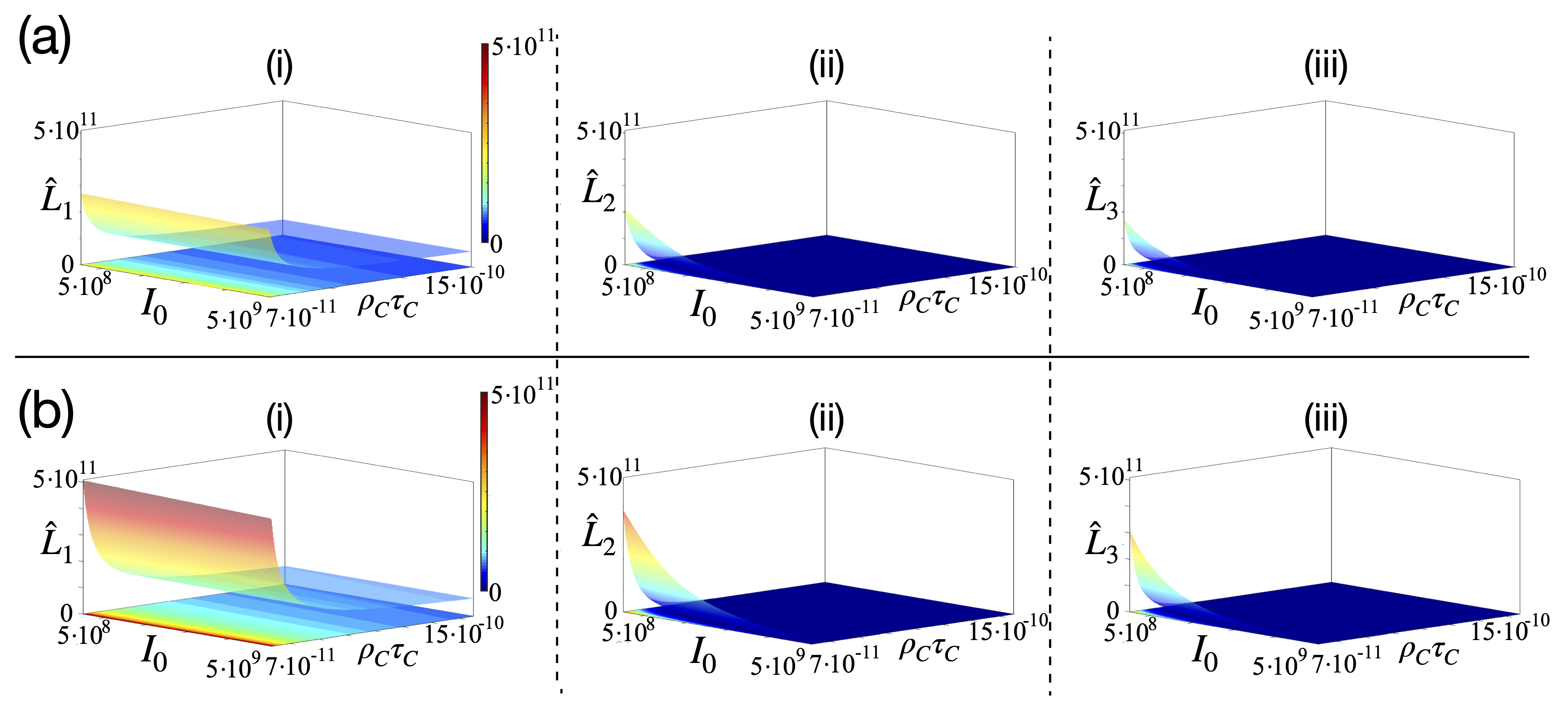}
\caption{Magnitude of the maxima of $L(t)$ as a function of $I_0$ and $\tau_C \rho_C$. The rest of parameter values are $\rho_L=0.2$ day$^{-1}\cdot$ cell$^{-1}$, $\alpha=10^{-11}$ day$^{-1}\cdot$ cell$^{-1}$, $\tau_I=4$ days and $\tau_B=45$ days. Initial conditions are $C_0=5\cdot 10^7$ cells, $L_0=5\cdot 10^{10}$ cells and $B_0=5\cdot 10^8$ cells. \textbf{(a)} Results for $\tau_C=14$ days.  \textbf{(b)} Results for $\tau_C=30$ days. Columns represent the first, second and third maxima. The color and coordinate $z$ represent the value of $L$ at the corresponding peak. \label{sura02-3}}
\end{figure*}

\begin{figure*}[ht!]
\centering
\includegraphics[width=0.85\textwidth]{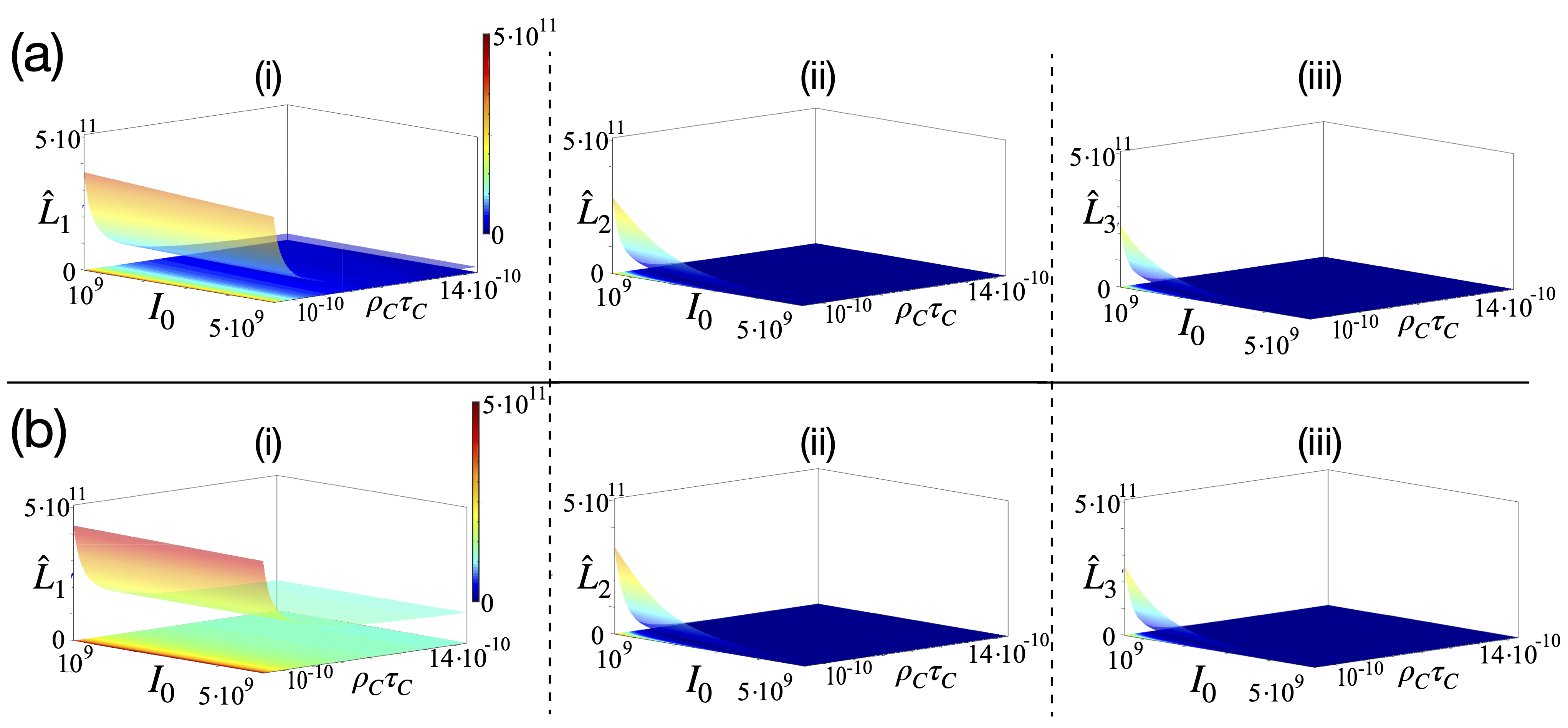}
\caption{Magnitude of the maxima of $L(t)$ as a function of $I_0$ and $\tau_C \rho_C$. The rest of parameter values are $\rho_L=0.2$ day$^{-1}\cdot$ cell$^{-1}$, $\alpha=10^{-11}$ day$^{-1}\cdot$ cell$^{-1}$, $\tau_I=4$ days and $\tau_B=45$ days. In this case, both rows have $\tau_C=22$ days and initial conditions $C_0=5\cdot 10^7$ cells and $B_0=5\cdot 10^8$ cells. \textbf{(a)} Results for $L_0= 10^{10}$ cells. \textbf{(b)} Results for $L_0= 10^{11}$ cells. Columns represent the first, second and third maxima. The color and coordinate $z$ represent the value of $L$ at the corresponding peak.}
\label{surat01}
\end{figure*}

\section{Analysis of successive maxima reached by the tumor}

\begin{figure*}[ht!]
\centering
\includegraphics[width=0.85\textwidth]{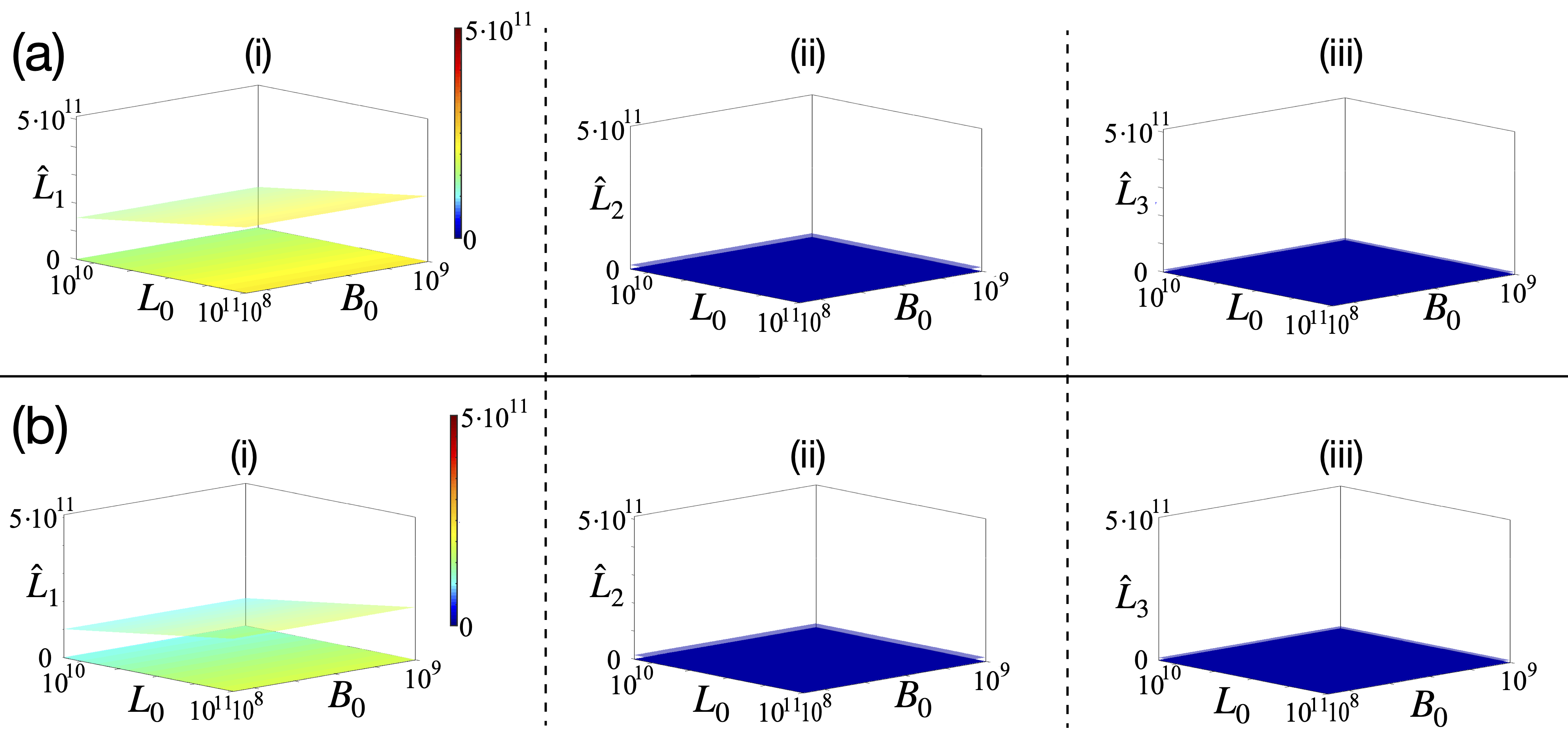}
\caption{Magnitude of the maxima of $L(t)$ as a function of initial conditions $L_0$ and $B_0$. Parameter values are $\rho_C=10^{-11}$ day$^{-1}\cdot$ cell$^{-1}$, $\tau_C=22$ days, $I_0=10^9$ cells, $\rho_L=0.2$ day$^{-1}\cdot$ cell$^{-1}$, $\alpha=10^{-11}$ day$^{-1}\cdot$ cell$^{-1}$, $\tau_I=4$ days and $\tau_B=45$ days. \textbf{(a)} Results for $C_0=10^7$ cells. \textbf{(b)} Results for $C_0=10^8$ cells. Columns represent the first, second and third maxima. The color and coordinate $z$ represent the value of $L$ at the corresponding peak.}
\label{cia01-3}
\end{figure*}

\begin{figure*}[ht!]
\centering
\includegraphics[width=0.85\textwidth]{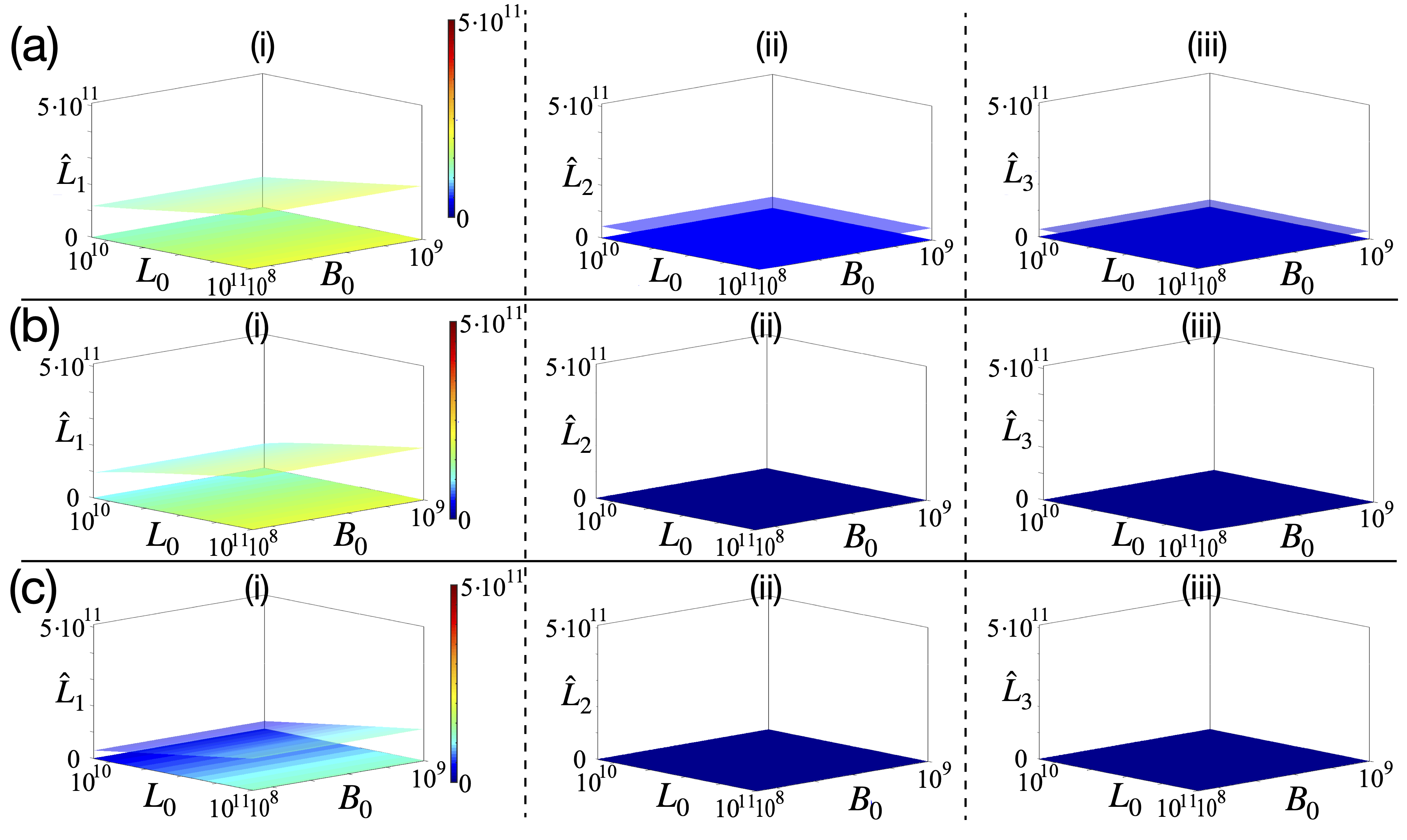}
\caption{Magnitude of the maxima of $L(t)$ as a function of initial conditions $L_0$ and $B_0$. Parameter values are $\tau_C=22$ days, $\rho_L=0.2$ day$^{-1}\cdot$ cell$^{-1}$, $\alpha=10^{-11}$ day$^{-1}\cdot$ cell$^{-1}$, $\tau_I=4$ days and $\tau_B=45$ days. $C_0=5\cdot 10^7$ cells in all cases. \textbf{(a)} Results for $\rho_C=10^{-11}$ day$^{-1}\cdot$ cell$^{-1}$ and $I_0=5\cdot 10^8$ cells. \textbf{(b)} Results for $\rho_C=10^{-11}$ day$^{-1}\cdot$ cell$^{-1}$ and $I_0=5\cdot 10^9$ cells. \textbf{(c)} Results for $\rho_C=5\cdot 10^{-11}$ day$^{-1}\cdot$ cell$^{-1}$, $I_0=5\cdot 10^8$ cells. Columns represent the first, second and third maxima. The color and coordinate $z$ represent the value of $L$ at the corresponding peak.}
\label{cia02b03c03}
\end{figure*}

\begin{figure*}[ht!]
\centering
\includegraphics[width=0.85\textwidth]{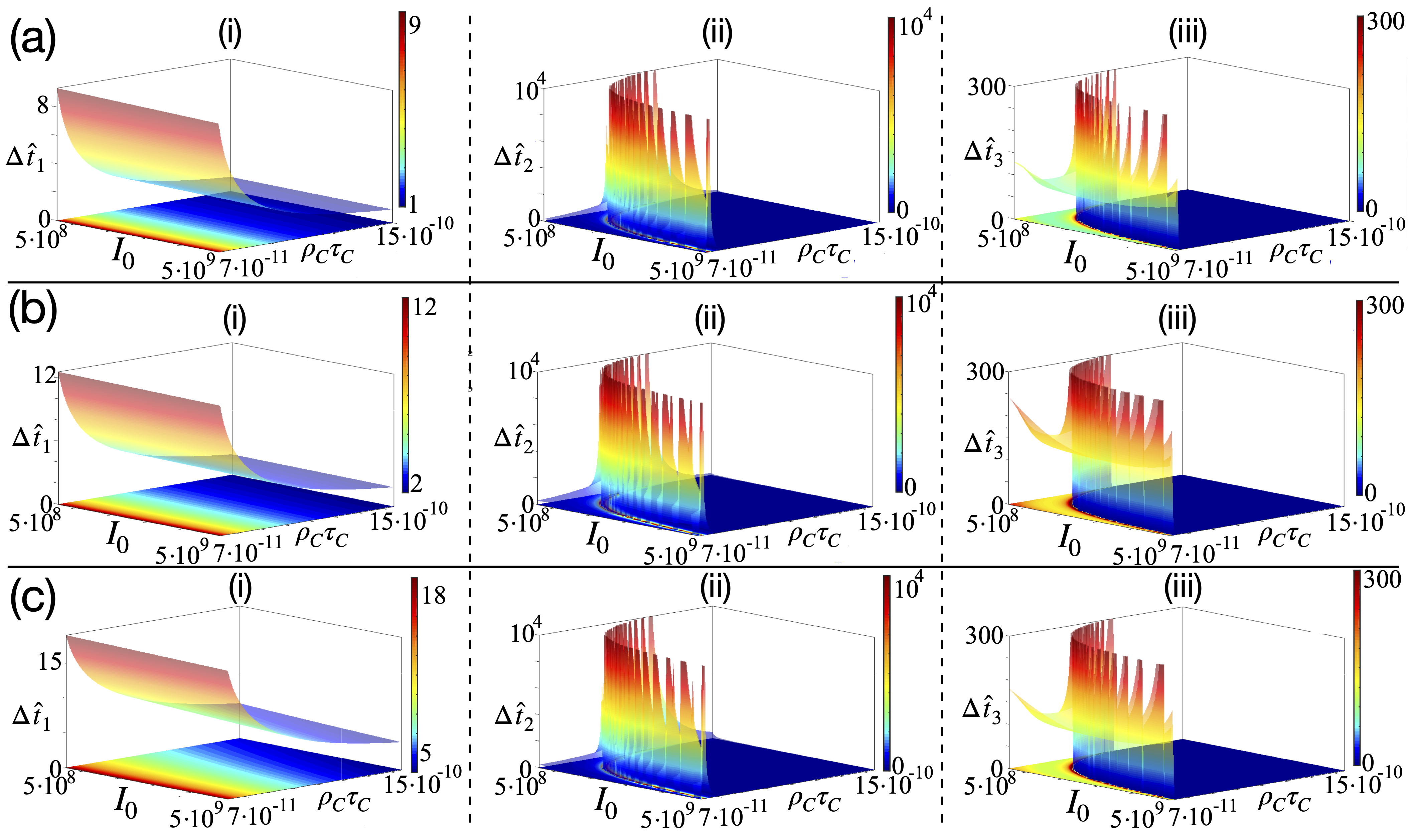}
\caption{Elapsed times between local maxima of $L(t)$ as a function of $I_0$ and $\rho_C \tau_C$. The rest of parameter values are $\rho_L=0.2$ day$^{-1}\cdot$ cell$^{-1}$, $\alpha=10^{-11}$ day$^{-1}\cdot$ cell$^{-1}$, $\tau_I=4$ days and $\tau_B=45$ days. Initial conditions are $C_0=5\cdot 10^7$ cells and $B_0=5\cdot 10^8$ cells. \textbf{(a)} Results for $\tau_C=14$ days and $L_0=5\cdot 10^{10}$ cells. \textbf{(b)} Results for $\tau_C=30$ days and $L_0=5\cdot 10^{10}$ cells. \textbf{(c)} Results for $\tau_C=22$ days and $L_0=10^{10}$ cells. In the first column the time it takes to reach the first peak from the start. In the second and third columns the time of the next peak from the previous one. These time values are represented (in days) both with the color scale and with the variable $z$.}
\label{tma02-3t01}
\end{figure*}

\begin{figure*}[ht!]
\centering
\includegraphics[width=0.85\textwidth]{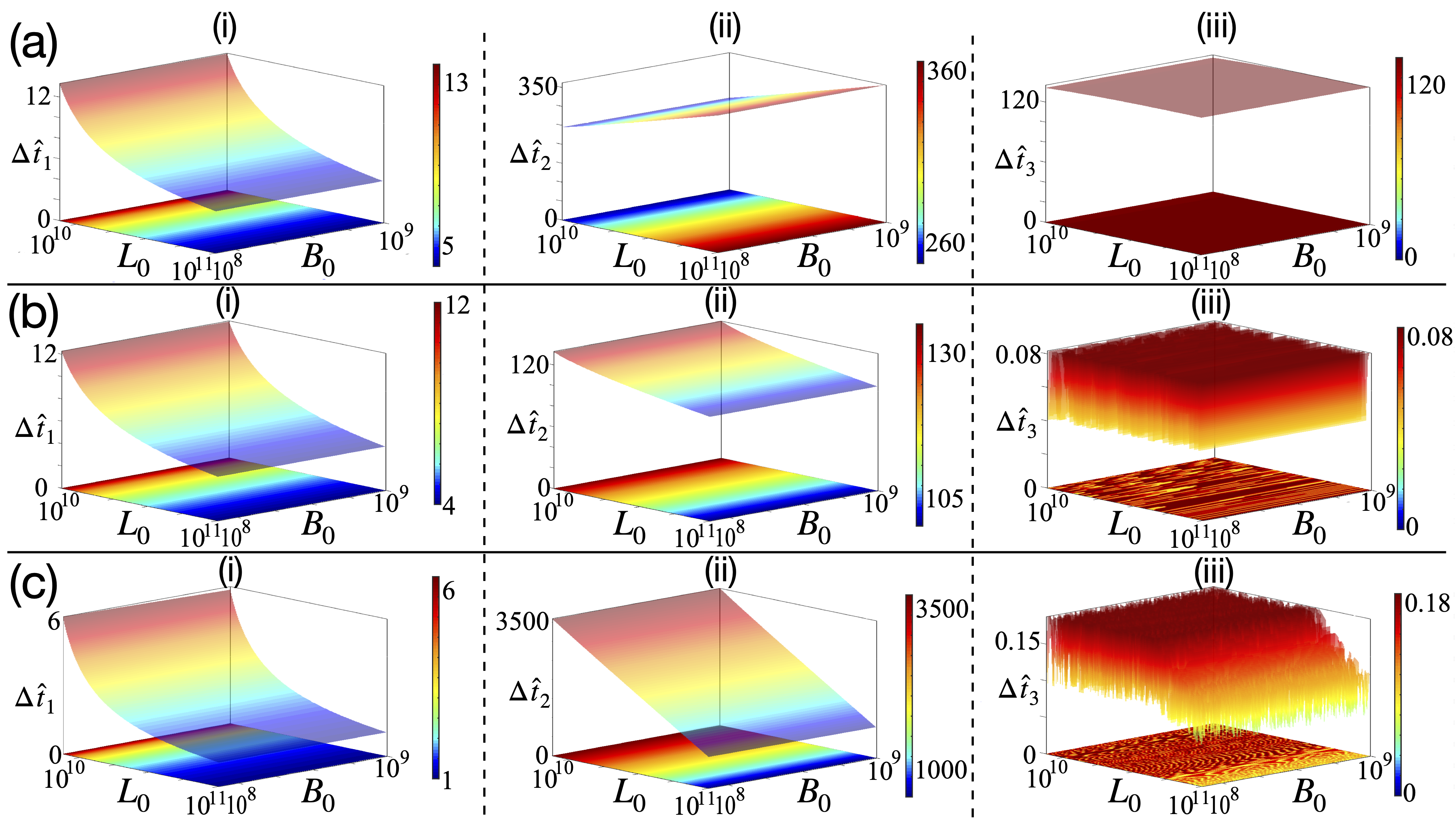}
\caption{Elapsed times between local maxima of $L(t)$ as a function of $L_0$ and $B_0$. The rest of parameter values are $\rho_L=0.2$ day$^{-1}\cdot$ cell$^{-1}$, $\alpha=10^{-11}$ day$^{-1}\cdot$ cell$^{-1}$, $\tau_I=4$ days and $\tau_B=45$ days. Here we have $\tau_C=22$ days and $C_0=5\cdot 10^7$ cells in all cases. \textbf{(a)} Results for $\rho_C=10^{-11}$ day$^{-1}\cdot$ cell$^{-1}$ and $I_0=5\cdot 10^8$ cells. \textbf{(b)} Results for $\rho_C=10^{-11}$ day$^{-1}\cdot$ cell$^{-1}$ and $I_0=5\cdot 10^9$ cells. \textbf{(c)} Results for $\rho_C=5\cdot 10^{-11}$ day$^{-1}\cdot$ cell$^{-1}$ and $I_0=5\cdot 10 ^8$ cells. In the first column the time it takes to reach the first peak from the start. In the second and third columns the time of the next peak from the previous one. These temporary values are represented (in days) both with the color scale and with the variable $z$.}
\label{tmcia02b03c03}
\end{figure*}

In this section we study the succession of maxima reached by the tumor, especially in the oscillatory regime described above. We display the results by means of three-dimensional plots, where the $z$ coordinate represents either the value of the peak (Figures \ref{sura02-3}-\ref{cia02b03c03}) or the time between successive peaks (Figures \ref{tma02-3t01}-\ref{tmcia02b03c03}). These values are also colormapped on the surface and projected on the $z=0$ plane. The $x$ and $y$ coordinates represent parameters or initial conditions that are varied in order to explore this property of the model.

Figure \ref{sura02-3} represents two different situations for the two extreme values of $\tau_C$. The free parameters are the same as those in Figure~\ref{equilibrios} to be able to relate the formulas obtained previously with the values of said maxima. In the top row, $\tau_C=14$ is taken. In the bottom row $\tau_C=30$. In both cases a similar behavior is observed. Although the actual values are somewhat different, especially for small values of $\tau_C\rho_C$, the overall behavior does not seem to be greatly affected by the value of $\tau_C$. The first column shows the first local maximum of $L$. This is not conditioned by the zones shown in Figure \ref{equilibrios} depending on the stability of the different equilibria. In fact, starting from a certain value $\tau_C\rho_C$ that is not very large, the variation is very little. Observe that in this case the influence of $I_0$ is practically non-existent. In the other two columns, in which the second and third local maximum are represented, respectively, the situation changes radically since now the influence of the stable equilibrium in the different parametric regions considerably determines the value of the maximum. In Figure \ref{surat01} we check the influence of initial tumor burden $L_0$. The top row takes $L_0=10^{10}$ and the bottom row takes $L_0=10^{11}$, with $\tau_C=22$ in both cases. The same differences between the first local maximum and the following ones are observed.

In the following two figures, the initial conditions $L_0$ and $B_0$ are left as free parameters instead of $I_0$ and $\tau_C \rho_C$. In Figure \ref{cia01-3} the parameters are set to standard values (placed in $R_2$ region) and two situations are considered, $C_0=10^7$ in the top row and $C_0=10^8$ in the bottom row. In both cases the numerical results are very similar, showing the little influence of $C_0$. The first maximum in $L$ is of the order of $2.5\cdot 10^{11}$, the second maximum drops to about $1.5\cdot 10^{10}$, and the third drops more slowly (about $10^{10}$). On the other hand, it can also be observed that $B_0$ hardly modifies the results either, being $L_0$ the only initial condition that influences, in an almost linear way, the values of the maxima that are produced initially. In Figure \ref{cia02b03c03} the parameter $C_0$ is fixed, $C_0=5\cdot 10^7$, and three different situations are studied: in the first row, $\rho_C=10^{-11}$, $I_0=5\cdot 10^8$ (placed in $R_2$); in the second row $\rho_C=10^{-11}$, $I_0=5\cdot 10^9$ (in $R_4$ region) while in the third $\rho_C=5\cdot 10^{-11}$, $I_0=5\cdot 10^8$ (in $R_3$).  The first case shows similar results to the Figure~\ref{cia01-3}, although the maxima  are slightly larger (by a factor of about 3). This situation is predictable since in both cases the parametric conditions are in the $R_2$ region, although in the latter case we are further away from the $R_3$ region and therefore the maxima in $L$ are larger (see Figure \ref{fococur}). The first maximum of the second situation examined behaves similarly. However, since we are now in $R_4$, the following maxima are practically negligible (of the order of unity). Something similar happens in the third situation. Now the first maximum is already somewhat lower than in the previous cases and the following maxima, being in $R_3$, are negligible.

Finally, the last two figures show the time that elapses between a local maximum and the next (or from the start at the first maximum). In Figure \ref{tma02-3t01} the parameters $I_0$ and $\tau_C\rho_C$ are left free for three different cases: firstly, standard initial conditions and extreme values of $\tau_C$ (14 and 30, respectively); lastly $L_0=10^{10}$ and $\tau_C=22$. In all cases, it is observed that the first maximum can take values from a few days to half a month. To reach the second peak, by contrast, can take several years if we are close to the region in which the eigenvalues of $P_3$ become real. This fits with the fact that the imaginary part tends to 0 and therefore the focus period to infinity (see Figure \ref{distancia}). The third maximum reproduces this situation, but with a much narrower band. This situation corresponds to the one shown in Figure \ref{het2eq}, in which the orbit passes close to the almost heteroclinic cycle between $P_2$ and $P_3$, staying for a long time close to $P_2$. In Figure \ref{tmcia02b03c03} the initial conditions $L_0$ and $B_0$ are left free and the same cases of Figure \ref{cia02b03c03} are studied. Again, we observe the key value $L_0$ and the difference between the first case (in $R_2$) and the next two in the second and third maxima. While in the first case the time between the first and the second maximum increases with $L_0$, in the other two cases it decreases, and in the third case the magnitude is also about 30 times higher. Moreover, the difference is even more obvious in the third column. The third maxima detected in the second and third cases correspond to small oscillations in the evolution, so the time difference between the second and third maxima is very small and irregular (hence the scatter plot with values close to 0). 

\section{Sensitivity analysis}

\begin{figure*}
\centering
\includegraphics[scale=0.5]{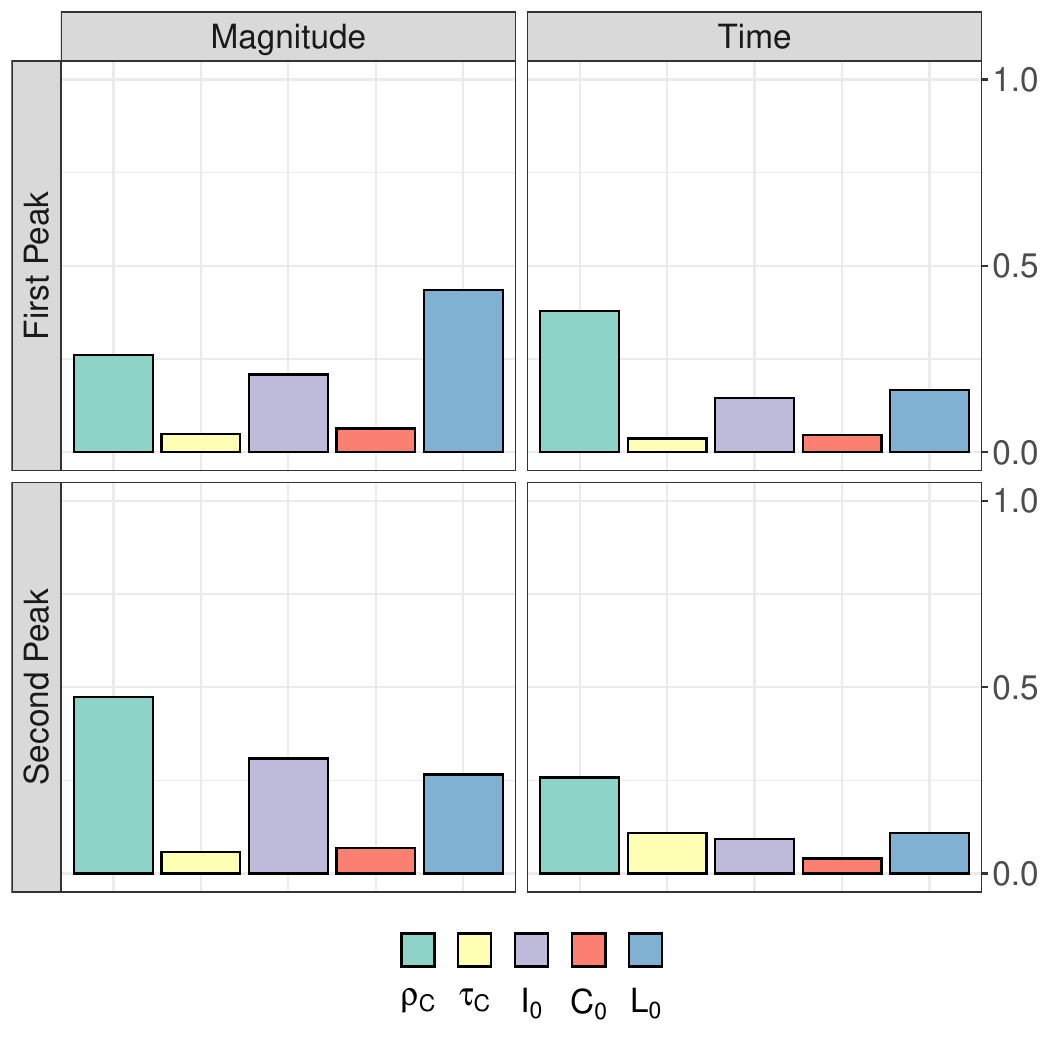}
\caption{Values of PAWN Sensitivity indices for magnitude and time of the first two maxima of $L(t)$. More information about their computation is given in Appendix~A.}
\label{SI}
\end{figure*}

In this section we perform sensitivity analysis to further support the results obtained in the previous one. This analysis measures simultaneously the impact of parameter variations on the outputs of interest, namely the magnitude of the first peaks of the leukemic compartments and the time it takes to reach them. We choose a moment-independent, global sensitivity analysis method called PAWN \cite{Pianosi15}. This  method measures the degree to which variations in one parameter induce changes in the distribution of the outputs of interest. The relative contributions of the parameters are what we call sensitivity index. They are displayed in Figure \ref{SI} for the magnitude of the first two peaks and the respective times it takes to reach such peaks. More details about the method can be found in Appendix~A.

We see that for the first peak the most important parameter is the initial tumor burden $L_0$, agreeing with what we found in the analysis of peaks previously performed. The influence of the bone marrow output $I_0$ is lower in this case. For the second peak, however, we stated that the region in the parameter space increases its importance, and we also see it here as $\rho_C$ becomes the parameter with the highest relative contribution, and $I_0$ surpasses $L_0$ in importance.

With respect to the timespan between maxima, we see that the CAR T-cell proliferation rate $\rho_C$ is the most relevant parameter, followed again by the bone marrow output $I_0$ and the initial tumor burden $L_0$. The conclusions are similar for both peaks, although the CAR T-cell mean life is more relevant in the second peak. The initial CAR T-cell number $C_0$ has low relative influence in general. These results confirm what was previously explored for discrete and non-simultaneous variations in the parameter ranges. 

\section{Discussion and conclusions}

In this work we have studied a modification of a previously published mathematical model of CAR T-cell therapy \cite{Leon21}. This mathematical model consisted of three ordinary differential equations for leukemic cells, CAR T-cells and B-cells that interacted with predator-prey-like dynamics: CAR T-cells recognize B-cells and leukemic cells, activate, proliferate, and remove these cells from the organism. The main intent of this paper was a more careful consideration of the input of B-cells coming from the bone marrow, which was included as a source term in the respective equation. This input plays an important role in this therapy for two main reasons. Since B-cells also share the CD19 marker with leukemic cells, the bone marrow input provides a continuous and sustained source of stimulation for anti-CD19 CAR T-cells. Early clinical trials already described the role of B-cells as an endogenous vaccine \cite{Heiblig16,Brudno16,Stein19}. On the other hand, CAR T-cell action induces permanent B-cell shortage or aplasia, a side-effect that is used as a surrogate marker of CAR T-cell persistence. In fact, the repopulation of the B-cell compartment is often followed by disease recurrence \cite{Davila16}. With this particular relevance in mind, we undertook stability analysis and bifurcation analysis to describe the dynamics of the system. We also carried out a description of the consecutive maxima reached by leukemic cells as well as the time at which these maxima are reached, in order to identify the most influential parameters of the model. This was completed by a sensitivity analysis by means of PAWN indices. 

After examining fundamental mathematical properties, our focus shifted to the biological implications of the mathematical model. In contrast to previous research on a simplified version of the model that identified unstable and stable equilibria involving CAR T-cells and leukemic cells in an oscillatory regime, our expanded model revealed an additional stable point where all three populations coexist. We observed that the continuous influx of B cells $I_0$ sustaining the pool of activated CAR T-cells could potentially lead to tumor eradication, depending on the values of $I_0$. Coexistence among the three populations occurs when $I_0<I^{\text{crit}}_0$, with a straightforward evolution toward focused dynamics on the 2D weak manifold. Conversely, if $I_0>I^{\text{crit}}_0$, the increased influx stimulates CAR T-cells enough to eliminate tumor cells completely. Thus, at $I_0=I^{\text{crit}}_0$, a transcritical bifurcation occurs, transforming the free-tumor equilibrium point $P_2$ into a stable point and the coexistence equilibrium point $P_3$ into an unstable point. These findings suggest the potential for pharmacological stimulation of stem cell asymmetric division and differentiation. An increase in $I_0$ could enhance the likelihood of complete cancer removal, aligning with hypotheses in other mathematical models of CAR T-cell therapy \cite{Martinez21}. A example of a therapeutic action that follow the same reasoning is the administration of T-cells modified to express CD19, which showed potential to improve the durability of CAR T-cell therapy in leukemia patients \cite{Annesley19}. 

We also found that the transition between $P_2$ and $P_3$ is mathematically described by an almost heteroclinic cycle between both points, with a slow movement in the vicinity of $P_2$ which then moves to the oscillatory regime around $P_3$. This behavior may explain the reason why after a long time where the illness is not present, it appears suddenly again (relapse). Heteroclinic-like cycles thus provide a successful theoretical explanation of this phenomenon, fulfilling one of the aims of the formal analysis of mathematical models of immunotherapy that we remarked in the introduction. 

Focusing on the oscillatory behavior of the system, we obtained an approximation for the time-dependency of both CAR T-cells and leukemic cells when sufficient time has passed. This allows for an analytic approximation of the asymptotic value of the successive maxima reached by the tumor, which was verified by means of numerical simulations. The magnitude of these maxima was later explored by means of biparametric interventions, first with $I_0$ and $\rho_C \tau_C $ and then with the initial conditions $L_0$ and $B_0$. We found that the influence of $I_0$ is little in this case, with $\rho_C$ begin a more relevant parameter. This means that B-cell input may determine the possibility of disease eradication, but it does not affect significantly the values reached by the tumor. Among the initial conditions, $L_0$, the initial tumor burden, was the only one to influence, almost linearly, the value of the maxima. This parameter was also the most relevant when studying the time that passes between the first, second and third maxima reached by the tumor. These results were later reinforced by PAWN sensitivity analysis, which incorporates simultaneous variations of several model parameters. The importance of tumor burden for therapy outcome was already suggested by clinical trials and later confirmed by real-word application of CAR T-cell therapy \cite{Larson21,Cappell23}.

Other important prognostic factors are the depth of the initial response and the administration of lymphodepleting chemotherapy prior to infusion. The goal of this therapy is to remove host T-cells, creating a favorable environment for CAR T-cells to expand. Its importance is gaining relevance recently and it could be included in the mathematical model by means of an additional compartment for the endogenous T-cell population. This has already been explored in other mathematical models \cite{Owens21,Kimmel21}. Other possible applications of the model studied here would be the exploration alternative therapeutic schemes or combination therapies, although the latter is more common in solid tumors. These tumors would require a different modelling approach since the spatial aspects become fundamental. In addition, we did not study the case of CD19$^-$ relapses. These happen when leukemic cells lose expression of CD19 and are therefore unaffected by CAR T-cell action. This phenomenon requires the consideration of evolutionary pressure of the treatment, for example by considering an additional CD19$^-$ compartment or partial differential equations with a variable for CD19 expression. Both have also been explored in other mathematical works \cite{Liu22,Santurio24}. Finally, while the focus of this paper was the formal analysis of the mathematical model, the field is already transitioning from proof-of-concept research to data-driven research. This opens the door to further modifications in the modelling framework and to studies of a more applied kind, which could help in verifying the results put forward by our analysis. 

In conclusion, we have studied a mathematical model describing the response of acute lymphoblastic leukemias to the administration of CAR T-cells, with a particular focus on the role of B-cell input from the bone marrow. Our study highlights the possibility of tumor eradication following stimulation of this bone marrow influx. It also provides a theoretical explanation of some observations reported in different clinical trials, such as the conceptualization of relapses as an heteroclinic connection between two equilibria. We hope this study stimulates the formal analysis of immunotherapy models and opens new research avenues in the area of CAR T-cell therapy.

\section*{Acknowledgments}
This work was partially supported by projects PID2022-142341OB-I00 and PID2022-140451OA-I00
funded by Ministerio de Ciencia e Innovación/Agencia Estatal de investigación (doi:10.13039/501100011033) and by ERDF A way of making Europe; by projects SBPLY/23/180225/000041 and SBPLY/21/180501/000145, funded by Junta de Comunidades de Castilla-La Mancha, Spain (and European Regional Development Fund (FEDER, EU)), and by University of Castilla-La Mancha / FEDER (Applied Research Projects) under grant 2022-GRIN-34405. AM has also been supported by Asociacion Pablo Ugarte. RB and SS have been supported by project PID2021-122961NB-I00 funded by Ministerio de Ciencia e Innovación/Agencia Estatal de investigación (doi:10.13039/501100011033) and by ERDF A way of making Europe, and by Diputaci\'on General de Arag\'on (grant numbers E24-23R and LMP94-21). 

\section*{Author Declarations}

\subsection*{Conflict of Interest}

The authors declare no competing interests.

\subsection*{Author Contributions}

\textbf{Sergio Serrano:} Formal analysis (lead), Software (lead), Visualization, Writing - original draft, Writing - review \& editing.
\textbf{Roberto Barrio:} Formal analysis (lead), Software (lead), Visualization, Writing - original draft, Writing - review \& editing.
\textbf{Álvaro Martínez-Rubio:} Conceptualization, Formal analysis, Methodology, Software, Visualization, Writing - original draft, Writing - review \& editing.
\textbf{Juan Belmonte-Beitia:} Methodology, Writing - original draft, Writing - review \& editing.
\textbf{Víctor M. Pérez-García:} Conceptualization, Methodology, Writing - review and \& editing.

\section*{Data Availability}

Data sharing is not applicable to this article as no new data were created or analyzed in this study

\section*{Appendix A. PAWN Sensitivity Analysis}

\begin{figure*}
\centering
\includegraphics[width=0.9\textwidth]{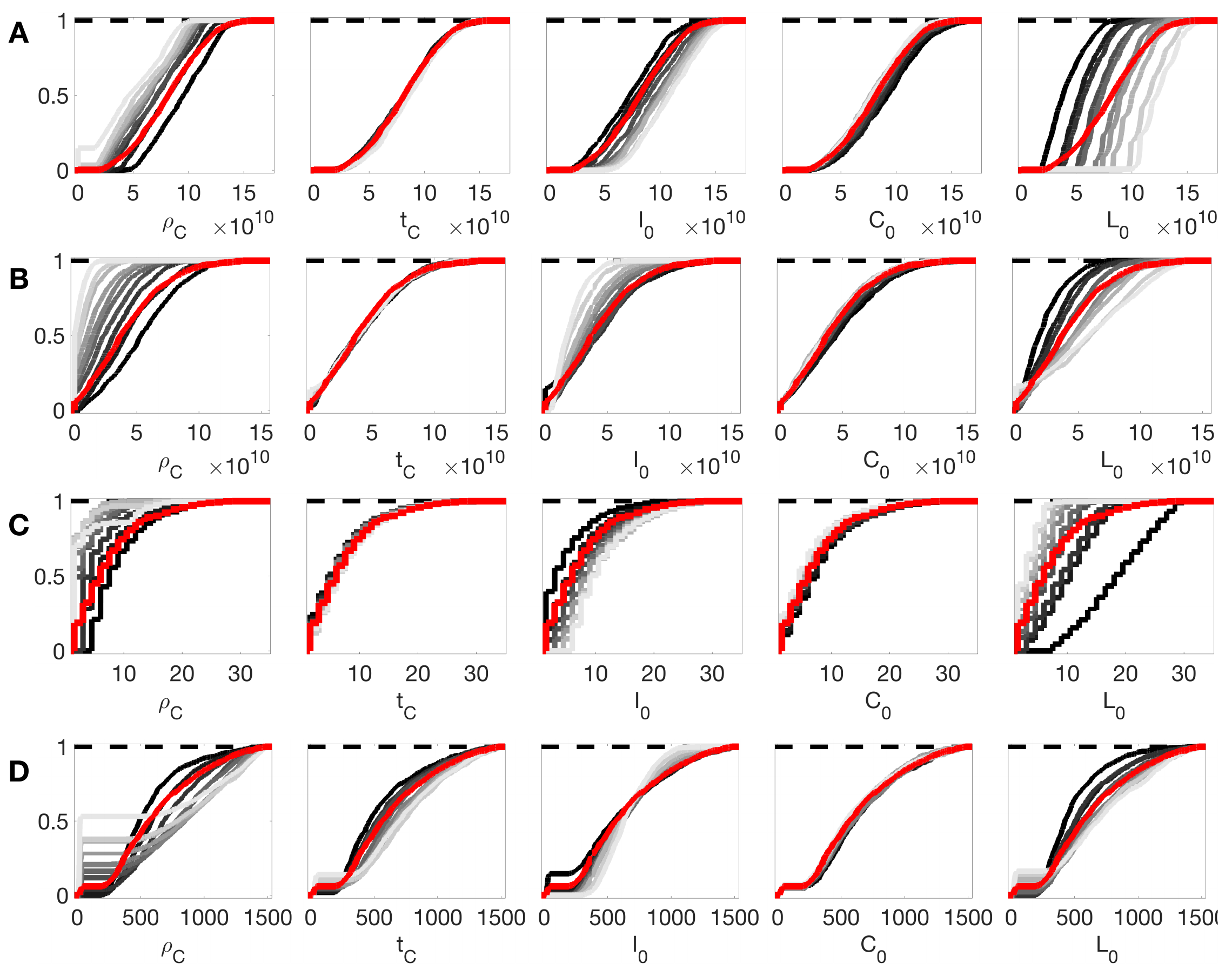}
\caption{PAWN method output distributions for the four outputs of interest: \textbf{(A)} Magnitude of first peak, \textbf{(B)} Magnitude of second peak, \textbf{(C)} Time to first peak, \textbf{(D)} Time to second peak. In red, unconditional distribution. In gray, conditional distributions for each conditioning value (colormap represents from lower (black) to higher value (white).}
\label{figS1}
\end{figure*}

There are many methods available for determining the propagation of uncertainty in a model. The most popular ones seem to be the variance-based approaches, which however have one limitation: They lead to contradictory results if the output distribution is highly skewed or multi-modal. This is our case, as can be seen in Fig \ref{figS1} where we show, in red, the unconditional distribution of the four outputs of interest. For this reason we resort to a recently developed moment-independent method \cite{Pianosi15}. The algorithm is available as an open source MATLAB code \cite{Pianosi15b}. 

\begin{figure*}
\centering
\includegraphics[width=0.9\textwidth]{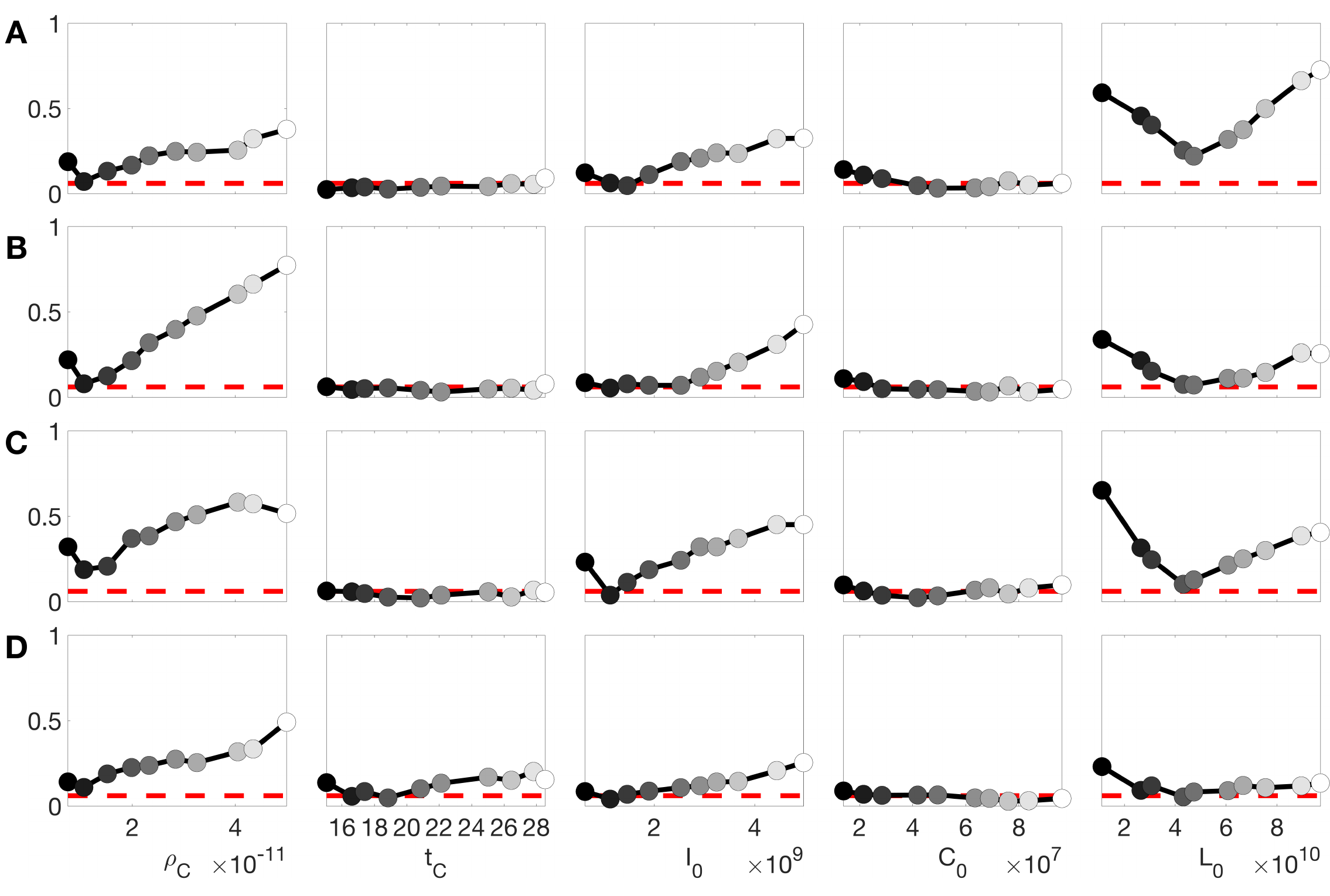}
\caption{Kolmogorov-Smirnov statistic for the conditioning points of each parameter of interest. \textbf{(A)} Magnitude of first peak, \textbf{(B)} Magnitude of second peak, \textbf{(C)} Time to first peak, \textbf{(D)} Time to second peak. The dashed red line represents the critical value of the KS statistic at confidence level of 0.05. The sensitivity index displayed in Figure \ref{SI} is the median of these values. The marker colormap is aligned with that of figure \ref{figS1}: Darker colors represent lower values of the parameter.}
\label{figS2}
\end{figure*}

The method works as follows. We first determine the output of interest. Then we select the parameters whose influence we aim to assess. In our case, these are CAR T-cells parameters $\rho_C$, $\tau_C$, the bone marrow output $I_0$ and the initial conditions $C_0$ and $L_0$. We specify a range of variation for each of them, which we take to be those used in Table~\ref{tabla}. The rest of parameters take the fixed values of Table~\ref{tabla}. We then compute the unconditional output distribution, by taking random points in the space of parameter values defined before. These are shown in red for each output in Figure~\ref{figS1}. We then compute conditional distributions by selecting random points while keeping one of the parameters constant, for different fixed values (conditioning points). These conditional distributions are shown in black in Figure~\ref{figS1}, one per each conditioning point. The Kolmogorov-Smirnov statistic is used to measure the difference between these distributions, and their values are shown in Figure~\ref{figS2}, again for each conditioning point. The median of these values is finally taken to be the sensitivity index, displayed as a relative value in Figure~\ref{SI}.
		
\section*{References}
\bibliography{bibliography}

\end{document}